\documentclass[journal,10pt]{IEEEtran}
\usepackage{mydef}

\makeatletter
\def\hlinewd#1{%
  \noalign{\ifnum0=`}\fi\hrule \@height #1 \futurelet
   \reserved@a\@xhline}
\makeatother

\pdfoutput=1
\begin{document}
\title{Distributed Coupled Multi-Agent Stochastic Optimization}

\author{
Sulaiman A. Alghunaim and     Ali H.~Sayed, ~\IEEEmembership{Fellow,~IEEE}
\thanks{A short preliminary conference version appears in \cite{alghunaim2018icassp}.}
\thanks{ S. A. Alghunaim is with the Department of Electrical Engineering, University of California, Los Angeles, CA 90095, USA (e-mail: salghunaim@ucla.edu). A. H. Sayed is with the Ecole Polytechnique Federale de Lausanne
(EPFL), School of Engineering, CH-1015 Lausanne, Switzerland (e-mail:
ali.sayed@epfl.ch). This work was supported in part by NSF grants CCF-
1524250 and ECCS-1407712. }}

\markboth{}%
{Shell \MakeLowercase{\textit{et al.}}: Bare Demo of IEEEtran.cls for Journals}

\maketitle

\begin{abstract}
This work develops effective  distributed strategies for the solution of constrained multi-agent stochastic optimization problems with coupled parameters across the agents. In this formulation, each agent is influenced by only a subset of the entries of a global parameter vector or model, and is subject to convex constraints that are only known locally. Problems of this type arise in several applications, most notably in disease propagation models, minimum-cost flow problems, distributed control formulations, and distributed power system monitoring. This work focuses on stochastic settings, where a stochastic risk function is associated with each agent and the objective is to seek the minimizer of the aggregate sum of all risks subject to a set of constraints. Agents are not aware of the statistical distribution of the data and, therefore, can only rely on stochastic approximations in their learning strategies. We derive an effective distributed learning strategy that is able to track drifts in the underlying parameter model. A detailed performance and stability analysis is carried out showing that the resulting coupled diffusion strategy converges at a linear rate to an $O(\mu)-$neighborhood of the true penalized optimizer.
\end{abstract}


\IEEEpeerreviewmaketitle

\section{Introduction}
\subsection{Motivation and Problem Setup}
In most multi-agent formulations of distributed optimization
problems, each agent generally has an individual cost function, $J_k(.)$, and the goal is to minimize the aggregate sum of the costs subject to some constraints, namely,
\begin{align}
 \underset{w \in \real^M}{\text{min   }}& \quad
  \sum_{k=1}^N J_k(w), \ \text{s.t.} \quad
 w \in \mathbb{W}_1 \cap \cdots \cap \mathbb{W}_N \label{glob1}
\end{align}
where $\mathbb{W}_k$ denotes a convex constraint set at agent $k$ and $N$ is the number of agents.  The aggregate  cost in \eqref{glob1} has one independent variable, $w\in  \real^{M}$,
which all agents need to agree upon \cite{6,34,35,36,56,23,24,25,29,33,5,7,28,8,27,26,37,38,39}. However, there exist many scenarios where
the global cost function may involve multiple independent
variables, and, moreover, each local cost may be a function
of only a subset of these variables. This situation motivates us to study in this work a broader
problem, where each local cost contains multiple variables that
get to be chosen by the network cooperatively. Examples of applications where this general scenario arises include  web page categorization \cite{14}, web-search ranking \cite{15}, disease progression modeling \cite{16}, distributed unmixing of hyperspectral data \cite{17,18}, minimum-cost flow problems \cite{19}, distributed model predictive control in smart energy systems \cite{11}, remote monitoring of physical phenomena involving discretization of spatial differential equations \cite{20}, distributed wireless acoustic sensor networks \cite{22}, distributed wireless localization \cite{31}, and distributed power systems monitoring \cite{21}.

 Thus, assume we have $L$ vector variables denoted by $\{w^1,\cdots,w^L\}$, where $w^\ell \in \mathbb{R}^{M_\ell}$. Let $w \triangleq {\rm col}\{
w^1, w^2, ... ,  w^L\} \in \mathbb{R}^{M}$ denote the $L \times 1$ block column vector formed by collecting all those variables. The partitioning is used to represent the possibility of multiple independent arguments for the cost functions. Without loss of generality, we assume that the variables $\{w^\ell\}$ are distinct in that they do not share common entries. 
 \begin{figure}[H]
\centering
\includegraphics[scale=0.5]{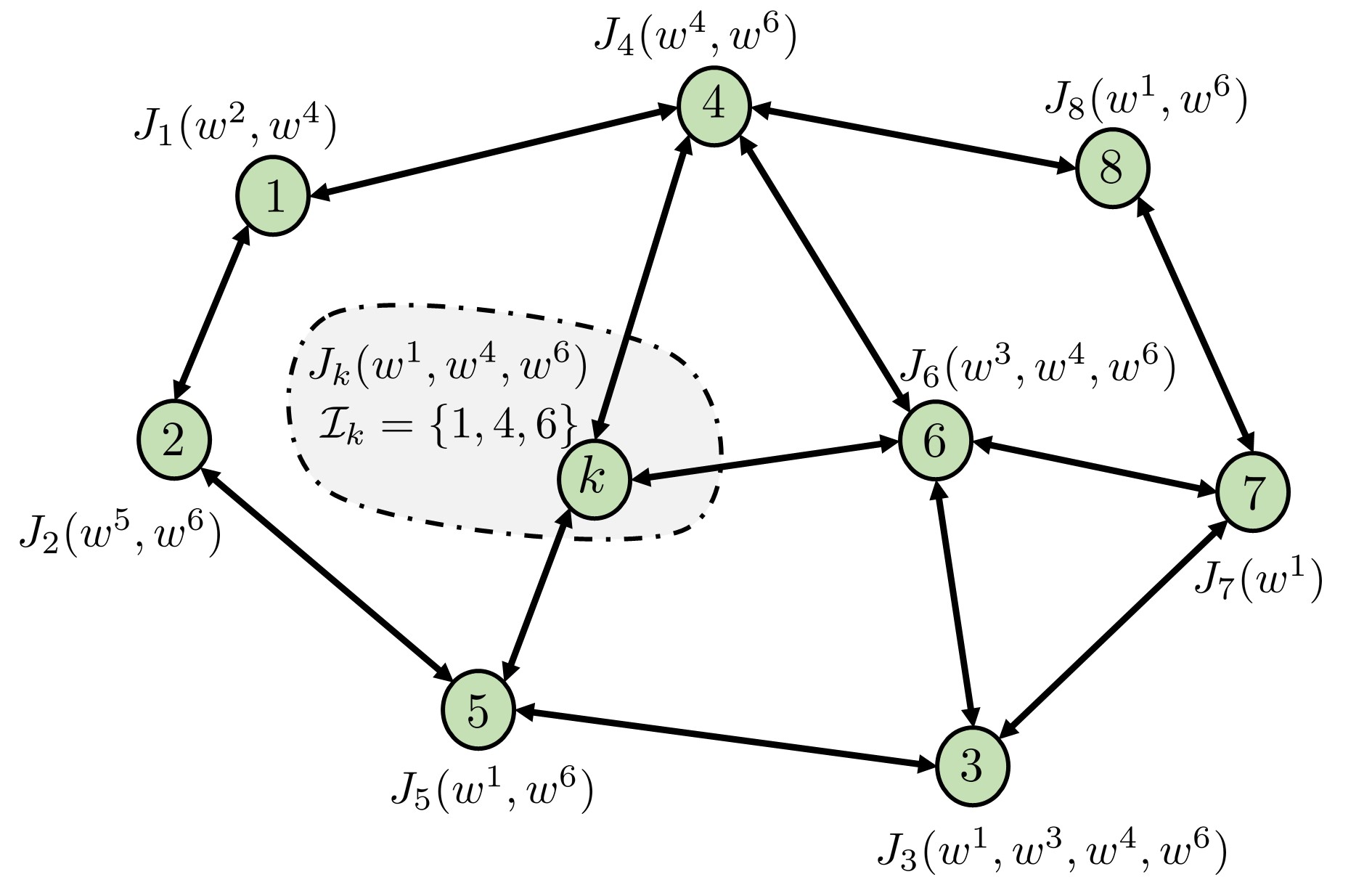}
\caption{A connected network of agents where different agents generally
depend on different subsets of parameter vectors. For this example, we have  $w=[w^1,w^2,w^3,w^4,w^5,w^6]$.}
\label{fig:network}
\end{figure} 
\noindent Let $\mathcal{I}_k$ denote the set of variable indices that affect the cost of agent $k$ -- Figure \ref{fig:network} illustrates this situation for a simple network. If we let $w_{k}$ denote the components of $w$ that affect this same agent:
\eq{
 w_k \triangleq {\rm col}\{w^{\ell}\}_{\ell \in \mathcal{I}_k} \in \mathbb{R}^{Q_k}, \quad Q_k \triangleq \sum_{\ell \in \cI_k} M_\ell.
\label{w_I_k}} 
Then we are interested in determining the solution of the following optimization problem:
\begin{align}
 \underset{w}{\text{min   }}& \ \
 J^{\rm glob}(w) \triangleq \sum_{k=1}^N J_k(w_k) \label{glob2}, \\ 
 \text{subject to    }& \ \
 w \in \mathbb{W}_1 \cap \cdots \cap \mathbb{W}_N \nonumber
\end{align}
The constraints set $\mathbb{W}_k$ is generally described by equality and inequality conditions of the form:
\eq{
\mathbb{W}_{k}=
\left\{
w:
\begin{array}{c}
 h_{k,u}(w_k)=0 ,\hspace{2 mm} u=1,.....,U_k\\
 g_{k,v}(w_k)\leq 0 , \hspace{2 mm} v=1,.....,V_k
\end{array}
\right.
}
where the functions $\{h_{k,u}(\cdot)\}$ are affine and $\{g_{k,v}(\cdot)\}$ are convex. 

 We note that algorithms that solve \eqref{glob1} can be used
to solve \eqref{glob2}. For example, this can be achieved by extending each local variable $w_k$ into the
longer global variable $w$. However, this solution method would require unnecessary
communications and memory allocation. This is because in
\eqref{glob2} each local function contains only a subset of the global
variable $w$. Therefore, solving \eqref{glob2} directly and more effectively is important for
large scale networks. Conversely, we also note that algorithms that solve
\eqref{glob2} are more general and can be used to solve \eqref{glob1}. To see
this, let $L = 1$ and ${\cal I}_k= \{L\}$, then problem \eqref{glob2} will
depend only on one variable $w=w^{L}$. In this case, the  cost function
becomes $J^{\rm glob}(w)=\sum_{k=1}^N J_k(w)$, which is of the same exact form
as problem \eqref{glob1}.
\begin{assumption}(\textrm{\bf{Strongly convex aggregate cost}})\label{feasible assump}: Problem \eqref{glob2} is feasible and has a strongly-convex cost: 
\eq{
\left(\grad_w J^{\rm glob}(x)-\grad_w J^{\rm glob}(y)\right)\tran (x-y) \geq \nu \|x-y\|^2 \label{stronglyconvex_bound}
}
for some constant $\nu>0$.
 \qd
\end{assumption}
\noindent Strong convexity is often assumed in the literature and is not a limitation since regularization is normally used to avoid ill-conditioning and it helps ensure strong convexity.  Under this condition, a unique solution $w^o \in \mathbb{W}$ exists. We denote its block entries by
\eq{
w^o= {\rm col}\{w^{1,o}, \cdots ,w^{L,o} \} \define \argmin_{w \in \mathbb{W}_1 \cap \cdots \cap \mathbb{W}_N} J^{\rm glob}(w) \label{optimal-original}
}
\begin{assumption} \label{cost-assump}
(\textrm{\bf{Individual costs}}): It is assumed that each cost function, $J_k(w_k)$, is differentiable with a Lipschitz continuous gradient:
\eq{
\big\|\grad_{w_k}J_{k}(x_k)-\grad_{w_k}J_{k}(y_k)\big\| \leq \delta_k \|x_k-y_k\| \label{indv-cost-1}
}
 for some positive constant $\delta_k$.
 \qd
\end{assumption}  
\noindent 
\begin{remark}{\rm \textbf{(Stochastic costs):} \label{remark-special-case}
The individual risk functions in Problem \eqref{glob2} can be stochastic in nature, such as $J_k(w_k)=\Ex Q(w_k;\x_{k,i})$, where $\Ex$ denotes the expectation operator over the distribution of the data, $Q(\cdot)$ is some convex loss function in $w_k$, and $\x_{k,i}$ is generic notation for the data at agent $k$ at time  $i$. Possible choices for $Q(\cdot)$ are quadratic losses, logistic losses, exponential losses, etc. Generally, the distribution of the data will not be available so that the local risk functions, $\{J_k(w_k)\}$ are not known beforehand by the agents. Our technique for solving Problem \eqref{glob2} will rely on the use of stochastic approximations for the gradient vectors of these risk functions, such as using the gradient of loss functions directly at every iteration (see \eqref{gradient-model} further ahead):
\eq{
\widehat{\nabla_{w_k}} J_k(\bzeta_{k,i}) \define \grad_{w_k} Q(\bzeta_{k,i};\x_{k,i}) \label{gradient-Q}
}
where $\bzeta_{k,i}$ will denote an intermediate estimate for $w_k$ at time $i$. These approximations introduce gradient noise, which refers to the difference between the true gradient vector and the above approximation. The challenge later will be to establish convergence despite the presence of this persistent noise. A special case of Problem \eqref{glob2} was treated in \cite{54}, which did not include constraint sets $\{\mathbb{W}_k\}$ and  assumed perfect knowledge of $J_k(\cdot)$  (i.e., there was no gradient noise).}
\qd
\end{remark}
\begin{remark}{\rm \textbf{(Useful case):} \label{remark-useful-case}
We illustrate a special case of \eqref{glob2}, which is common in many applications. Consider the scenario where every agent wants to estimate its own variable $w^k$ and is coupled with every neighboring agent, i.e.,  $L=N$ and $\cI_k=\cN_k$ so that $w_k={\rm col}\{w^\ell\}_{\ell \in \cN_k}$, where $\cN_k$ denotes the neighborhood of agent $k$ (including agent $k$). To explicitly indicate that $w_k$ and the corresponding constraint depend exclusively on the neighborhood variables, we let
\eq{
w_{\cN_k}\define w_k , \ \mathbb{W}_{\cN_k} \define \mathbb{W}_{k} }
 Then, problem \eqref{glob2} becomes
\begin{align}
 \underset{w}{\text{min   }}& \quad
  \sum_{k=1}^N J_k(w_{\cN_k})\label{glob-application} 
 , \ \text{s.t.     } \quad
 w \in \mathbb{W}_{\cN_1} \cap \cdots \cap \mathbb{W}_{\cN_N}
\end{align}
Many important applications fit into problem \eqref{glob-application}. For, instance, such formulations arise in wireless localization where each agent aims to estimate its position based on distance measurements from its neighbors. Two other examples are distributed model predictive control \cite{52} and minimum cost flow problems \cite{19}.} \qd
\end{remark}
\noindent \begin{example}{\rm \textbf{(Power system state estimation):} \label{example-application}
 We describe one example in power system state estimation \cite{21}, which is a special case of formulation \eqref{glob2}. Thus, consider a system consisting of $N$ interconnected sub-systems with each sub-system consisting of some subset of buses (or edges). Let $w=[w^1,\cdots,w^L]$ denote the state of the system (e.g., voltages and currents across all buses). Suppose each subsystem collects measurements related to the voltages and currents across its local buses and voltages and currents across the interconnection between neighboring sub-systems (see Figure \ref{fig:power-system}). 
\begin{figure}[H]
\centering
\includegraphics[scale=0.6]{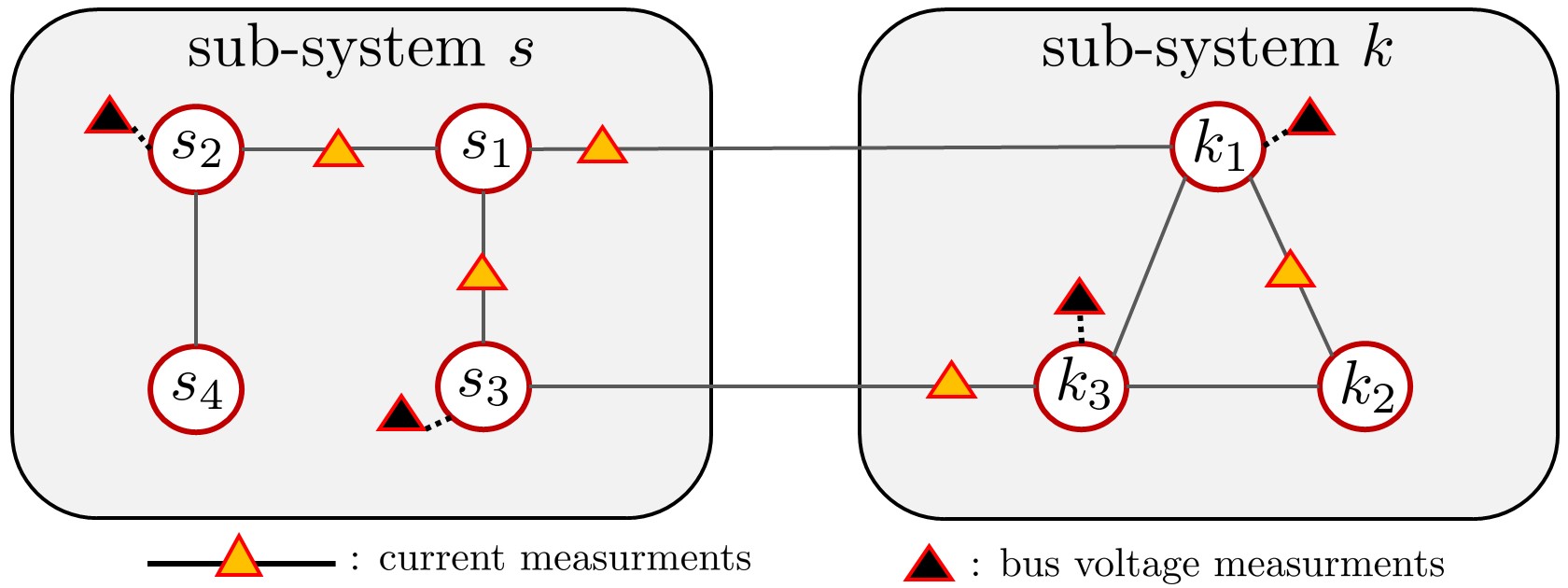}
\caption{Two neighboring sub-systems sharing states across their interconnection, i.e., buses $k_1$-$s_1$ and $k_3$-$s_3$.}
\label{fig:power-system}
\end{figure} 
 We let $w_k$ denote a vector that collects the states $\{w^\ell\}$ of system $k$ (i.e., voltages and currents across the buses of system $k$ and across the buses to its neighboring subsystems). Then, the goal of each subsystem is to estimate the states $w_k$ from $S_k$ observations:
\eq{
y_k=H_k w_k + v_k
}
where $H_k \in \real^{S_k \times Q_k}$ is the measurement matrix and $v_k \in \real^{S_k}$ is a zero-mean measurement noise with known covariance matrix. One way to estimate the $\{w_k\}$ is by solving the following problem:
\begin{align}
 \underset{w}{\text{min   }}& \quad
  \sum_{k=1}^N \|y_k-H_k w_k\|^2\label{power-system-estim} , 
 \ \text{s.t.    } \ \
 w \in \mathbb{W}_1 \cap \cdots \cap \mathbb{W}_N
\end{align}
where $\{\mathbb{W}_k\}$ are a convex sets that capture some prior information about about the $\{w_k\}$. Now, since neighboring agents measure some similar quantities across their interconnections, it holds that $w_k$ and $w_s$ are expected to partially overlap if $s \in \cN_k$ and hence this problem is a special case of \eqref{glob2}. 
}

\qd
\end{example}
\noindent \begin{example}{\rm \textbf{(Network flow optimization)\cite{19}:} We provide a second example that fits into problem \eqref{glob2}; a third example involving model-fitting when different subsets of data are dispersed over different agents is discussed in the simulations section. Consider a directed network with $N$ agents and $L$ links. Let $w^\ell$ denote the net flow in link $\ell$. Let $b_k$ denote an external supply (or demand) of flow for node $k$ such that $\sum_{k=1}^N b_k=0$ (i.e., flows entering the network match the flows leaving the network). Let $\cI_k$ denote the links connected to node $k$ so that $w_k={\rm col}\{w^\ell\}_{\ell \in \cI_k}$ denotes the vector of flows across the links connected to node $k$. Then, we can formulate the problem where the network of agents is interested in solving:
 \begin{align}
 \underset{w}{\text{min}} \
  \sum_{k=1}^N J_k(w_k) , \quad 
 \text{subject to     }& \quad
a_k\tran w_k =b_k, \quad \forall \ k \label{appl-flow} \\
& \quad   w \in \mathbb{W}_1 \cap \cdots \cap \mathbb{W}_N \nonumber
\end{align}
where $a_k\tran w_k =b_k$ is a flow conservation constraint such that $a_k$ is a vector with entries $+1$ and $-1$ with $+1$ at the position of entering flows and $-1$ at the position of leaving flows -- see Figure \ref{fig:flow-system}. Moreover, $J_k(\cdot)$ is some convex function and the $\{\mathbb{W}_k\}$ denote some convex constraints such as link capacity constraints -- see \cite{19}.
\begin{figure}[H]
\centering
\includegraphics[scale=0.4]{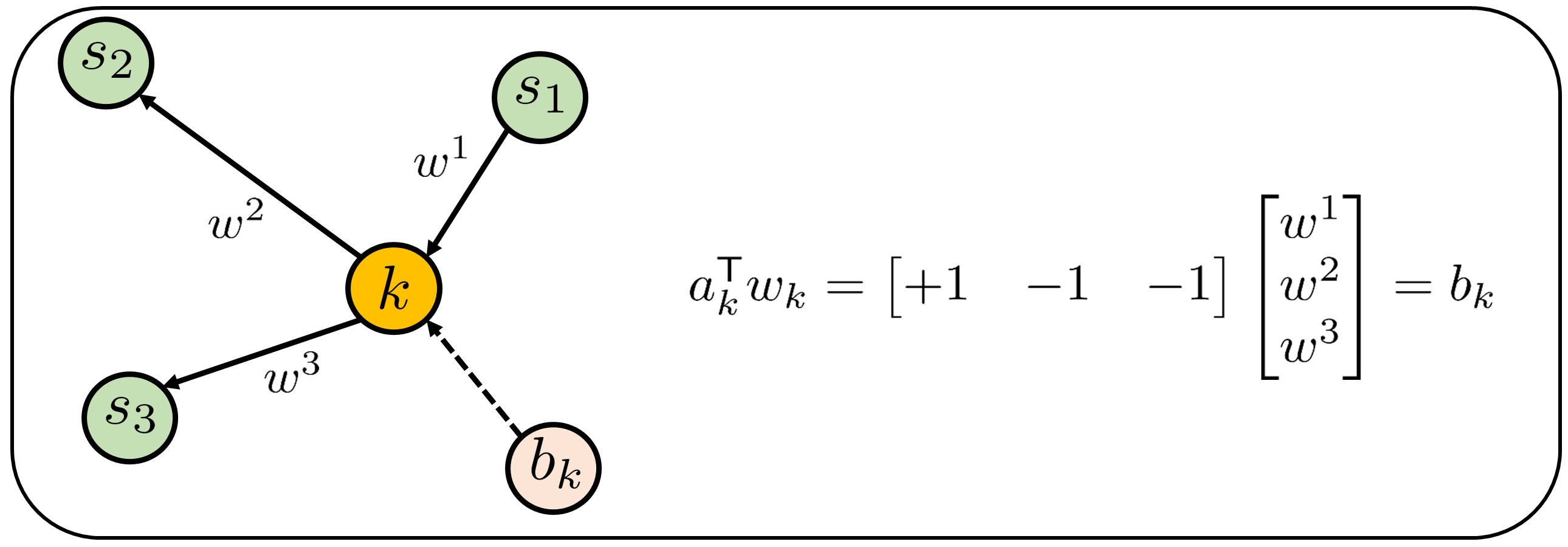}
\caption{An illustrative example to show how the constraint $a_k\tran w_k =b_k$ is formed.}
\label{fig:flow-system}
\end{figure} 
 Now we note that it is usually assumed that $b_k$ is measured exactly. However, in practice noise is usually present in $b_k$ and thus it can be modeled as $\b_k(i)=a_k\tran w_k+\v_k(i)$ where $\v_k(i)$ is some unknown measurement noise. If this is the case, then the constraint $b_k=a_k\tran w_k$ cannot be satisfied and one approach to address this situation is to employ a penalty method to solve instead \cite{1}:
 \begin{align}
 \underset{w}{\text{minimize   }}& \quad
  \sum_{k=1}^N J_k(w_k)+\eta \Ex\|a_k\tran w_k -\b_k(i)\|^2 \label{appl-flow-penalty} \\
 \text{subject to     }& \quad
 w \in \mathbb{W}_1 \cap \cdots \cap \mathbb{W}_N \nonumber
\end{align}
for some finite large penalty $\eta>0$.
}
\qd
\end{example}
\subsection{Contribution and Related Work}
There has been extensive work in the literature on solving problems of the type \eqref{glob1}, which would therefore be applicable to problem \eqref{glob2} albeit by going through the costly step of extending the local vectors, as explained before. For example, incremental strategies have been used in \cite{6,34,35,36}, consensus strategies in \cite{23,24,25,29,33}, diffusion strategies in \cite{5,7,28,8,27}, and ADMM strategies in \cite{26,37,38,39}. For sparse networks with a large number of parameters to estimate, it is much more efficient to devise distributed techniques that solve \eqref{glob2} {\em directly} rather than transform \eqref{glob2} into the form in \eqref{glob1} via vector extension. It is shown in the simulations in \cite{54} and \cite{12} that this extension technique not only increases complexity but it often degrades convergence performance as well, which we show analytically in this work.

Therefore, it is desirable to address the solution of problem \eqref{glob2} directly. Problems of this type have received less attention in the literature. For example, in deterministic formulations, ADMM techniques have been used to solve $\eqref{glob2}$ or its special case \eqref{glob-application} in \cite{21,10,12,32}.
In particular, the work \cite{21} applies an ADMM method to solve a distributed
power system state estimation problem of the form \eqref{glob2}, while the work \cite{12} solves \eqref{glob2} by employing an
extended ADMM method to reduce communications at the expense of some stronger assumptions. Likewise, in the model
predictive control literature, most of the methods used are
specific for the special case \eqref{glob-application} \cite{52}. For example, to solve \eqref{glob-application} in \cite{10} another
ADMM method is proposed, while \cite{32} uses an inexact
fast alternating minimization algorithm; this second method is
equivalent to an inexact accelerated proximal-gradient method applied to a dual problem \cite{48}. In all of these methods,
a second auxiliary (sub-minimization) problem needs to be solved at each iteration, which
requires an inner iteration unless a closed form solution exists. In the stochastic optimization literature, some special cases of \eqref{glob2} have also been considered. For example, the work \cite{4} focuses on multi-task unconstrained quadratic problems and employs game-theoretic techniques. In \cite{1,2} another quadratic problem is solved, where every agent has their own variable $w^k$ (i.e, $L=N$) and the agents are coupled through linear constraints with neighboring node's variables $\{w^\ell\}_{\ell \in \cN_k}$. Moreover, it is further assumed that the agents involved in a constraint are fully connected, i.e., they can communicate directly.

 Since in problem \eqref{glob2} different agents are influenced by different block vectors $w^\ell$, the network will be divided into overlapping clusters and each cluster $\ell$ will involve the agents that need to agree on $w^\ell$ -- see \eqref{cluster}. Similar clustering was used in \cite{12} for the deterministic problem where the ADMM method was employed with identical penalty factors across all clusters. In our previous work \cite{54}, we studied the same deterministic case but developed instead a first-order method for solving \eqref{glob2} without constraints by relying on the exact diffusion strategy from \cite{45,55}, which unlike ADMM does not require inner minimization steps.
 
 There are also works that deal with problem \eqref{glob1} where all agents need to agree on the same $w$, however to reduce communication, different agents transmit different blocks $\{w^\ell\}$ at each time instant \cite{arablouei2014distributed,notarnicola2017distributed}. In this case, each cluster involves different agents at each time instant and, over time, all agents will be involved in all clusters.  Note also that the group diffusion algorithm used in \cite{47} deals with the problem where each agent is interested in its own minimizer $
w^{\bullet}_k=\argmin_{w} J_k(w)$ and can solve the problem individually but cooperation is used since the estimation accuracy can be enhanced by cooperation if part of the minimizers $\{w^{\bullet}_k\}$ are common across neighbors. To take advantage of this overlap, each agent assigns different weights to different blocks  $w^\ell$.

In this work, we will solve problem \eqref{glob2} under a {\em stochastic} environment where agents do not necessary know the exact gradient information but are subject to noisy perturbations. We will also employ {\em constant} step-size learning in order to endow the resulting recursions with adaptation abilities to drift in the models. It was shown in \cite{40,41} that under such scenarios, diffusion strategies have superior performance than consensus strategies and primal-dual methods over adaptive networks. The superiority is due to some inherent asymmetries that exist in the updates of these latter methods, and which cause degradation in performance when the algorithms are required to learn continually from streaming data; the details and the origin of this behavior are explained in \cite{40,41}. Additionally, it is explained in \cite{45} that diffusion strategies can be motivated by optimizing penalized costs. For these reasons, we shall employ in this work penalized diffusion methods to solve problem \eqref{glob2}.  We prove that the algorithm converges linearly to an $O(\mu)$ neighborhood of the solution of a penalized cost, which can be made arbitrarily close to the solution of the original problem. One important conclusion from our analysis is to clarify the effect of the clustering step on the convergence rate. This was observed in \cite{12,54} through simulation and is explained analytically in this work.

While it is common in the literature to employ projection methods to solve constrained problems of the form \eqref{glob1}, by continuously projecting the iterated estimates onto the convex constraint sets, these methods nevertheless require the projection sets to be geometrically simple so that the projections can be computed efficiently. When this is not the case, penalty-based methods become more efficient. For example, in large-scale Markov decision problems (MDPs) exact dynamic programming techniques are not feasible since the computational complexity scales quadratically with the number of states. As such, in \cite{44}  a stochastic penalty based method is suggested to reduce the computational cost for large MDPs. Penalty-based methods are also attractive when the constraints are stochastic in nature, as happens in statistical estimation \cite{30} and stochastic minimum-cost flow problems \cite{1}. They are also useful when the constraints are soft in that they need not be satisfied exactly or when the constraints are used to encode prior information about the unconstrained optimal value.

\noindent {\bf Notation}.  We use boldface letters to denote random quantities and regular font to denote their realizations or deterministic variables. All vectors are column vectors unless otherwise stated. We use ${\rm col}\{x_j\}_{j=1}^{N}$ to denote a column vector formed by stacking $x_1, ... , x_N$ on top of each other, $\text{diag}\{x_j\}_{j=1}^{N}$ to denote a diagonal matrix consisting of diagonal entries $x_1, ... , x_N $, and $\text{blkdiag}\{X_j\}_{j=1}^{N}$ to denote a block diagonal matrix consisting of diagonal blocks $X_1, ... , X_N $. The notation $x \prec y$ ($x \preceq y$) means each entry of the vector $x$ is less than (or equal) to the corresponding vector $y$ entry. We use $O(\alpha)$ to indicate values of the order the scalar $\alpha$ (i.e., $O(\alpha)=c\alpha$ for some constant $c$ independent of $\alpha$). For any set $\cZ=\{z_1,z_2,\cdots,z_r\}$, where $z_1<z_2<\cdots<z_r$ are integers. We let $U=[G_{mn}]_{m,n \in \cZ}$ denote the $r \times r$ matrix with $(i,j)$-th entry equal to $G_{z_i z_j}$.
\section{Problem Formulation}
\subsection{Penalized Formulation}
In what follows we adjust the formulation from \cite{5} to problem \eqref{glob2}. One key difference in relation to what was studied in \cite{5} is that problem \eqref{glob2} allows individual agents to depend on {\em different} subsets of the parameter vector. That is, {\em coupling} now exists between the individual cost functions through the sharing of common sub-vectors. In contrast, in the formulation studied in \cite{5}, all agents were assumed to share the {\em same} parameter vector $w$ and were interested in reaching agreement about it. Here, instead, agents share {\em different} components of the larger vector $w$ and, through cooperation, they need to arrive at agreement. We shall develop a distributed scheme that enables this coupled learning objective and analyze its performance. In addition, the analysis in \cite{5} assumes each individual cost is twice-differentiable and strongly convex. Unlike \cite{5}, we do not impose any convexity assumption on the individual costs. We only require the aggregate cost to be strongly convex and Assumption \ref{cost-assump} to be satisfied. As a result the analysis becomes more challenging and is substantially different from the techniques used in \cite{5}.

One of the initial steps used in \cite{5} is to replace the constrained problem \eqref{glob2} by an unconstrained problem through the introduction of a penalty term; the purpose of this term is to penalize deviations from the constraints. Thus, we start by relaxing problem \eqref{glob2} by the following penalized form parametrized by a scalar $\eta > 0$:
\begin{align}
 \underset{w}{\text{minimize   }}& \quad
 J^{\rm glob}_{\eta}(w) \triangleq \sum_{k=1}^N J_{k,\eta}(w_k) \label{penalized_cost} 
\end{align}
where the individual costs on the right-hand side incorporate a penalty term, as follows:
\eq{
J_{k,\eta}(w_k) \triangleq  J_k(w_k)+\eta \hspace{1mm} p_k(w_k)}
with each penalty function given by
\eq{
 p_k(w_k)\triangleq \sum_{u=1}^{U_k} \delta^{\rm EP}\big(h_{k,u}(w_k)\big) +\sum_{v=1}^{V_k}\delta^{\rm IP}\big(g_{k,v}(w_k)\big) }
Here, the symbols $\delta^{\rm EP}(x)$ and $\delta^{\rm IP}(x)$ denote differentiable  convex functions chosen by the designer in order to  penalize the violation of the constraints, namely, they are chosen to satisfy the requirements:
 \eq{ \label{penalty-f}
\delta^{\rm EP}(x)=
\left\{\scalemath{0.9}{
\begin{array}{c}
 0 ,\hspace{5.5 mm} x=0\\
 >0 ,\hspace{2 mm} x\neq0
\end{array},}
\right.\quad 
 \delta^{\rm IP}(x)=
\left\{\scalemath{0.9}{
\begin{array}{c}
 0 ,\hspace{5 mm} x\leq0\\
 >0 ,\hspace{2 mm} \text{otherwise}
\end{array}}
\right.}
For example, the following two continuous, convex, non-decreasing, and twice differentiable choices that satisfy \eqref{penalty-f} are given in \cite{5}:
\eq{
\delta^{\rm EP}(x)= x^2, \quad
\delta^{\rm IP}(x)= \max\bigg(0, {x^3\over \sqrt{x^2+\rho^2}}\bigg)
}
for some parameter $\rho>0$. 
\begin{assumption} \label{penalty-assump}
(\textrm{\bf{Penalty functions}}): The penalty function $p_k(w_k)$ is convex and differentiable with a Lipschitz continuous gradient:
\eq{
\big\|\grad_{w_k}p_{k}(x_k)-\grad_{w_k}p_{k}(y_k)\big\| \leq \delta_{p,k} \|x_k-y_k\|  \label{penalty-cost-1}
}
for some positive scalar $\delta_{p,k}$.\qd
\end{assumption} 
Since $J^{\rm glob}(w)$ is strongly convex and each $p_k(w_k)$ is convex, the cost $J^{\rm glob}_{\eta}(w)$ is also strongly convex and, thus, a unique solution exists for problem \eqref{penalized_cost}, which we denote by:
\eq{
w^\star= {\rm col}\{w^{1,\star}, \cdots ,w^{L,\star} \} \define \argmin_{w^1,\cdots ,w^L} J^{\rm glob}_{\eta}(w) \label{optimal-penalized}
}
\subsection{Centralized Solution}
We first show how to solve problem \eqref{penalized_cost} in a centralized manner, assuming the presence of a central processor with knowledge of each individual cost $J_{k,\eta}(.)$. We  define:
\eq{
p^{\rm glob}(w) \define \sum_{k=1}^N p_k(w_k)
}
 Then, the gradient vector of \eqref{penalized_cost}  relative to $w = {\rm col}\{w^1, w^2, ... ,  w^L\}$ is given by:
\eq{
\grad_w J_{\eta}^{\rm glob}(w)&=\grad_w J^{\rm glob}(w) +\eta \grad_w p^{\rm glob}(w) \nonumber \\
&= \begin{bmatrix}
 \sum\limits_{k=1}^{N} \grad_{w^1}J_{k}(w_k) \\
\vdots \\
 \sum\limits_{k=1}^{N} \grad_{w^L}J_{k}(w_k)
\end{bmatrix} + \eta \begin{bmatrix}
 \sum\limits_{k=1}^{N} \grad_{w^1}p_{k}(w_k) \\
\vdots \\
 \sum\limits_{k=1}^{N} \grad_{w^L}p_{k}(w_k)
\end{bmatrix} \label{centralized-gradient}
}
Using a gradient descent algorithm, we can solve \eqref{penalized_cost} iteratively as follows:
\eq{
w_i=w_{i-1} - \mu D \grad_w J_{\eta}^{\rm glob}(w_{i-1}), \quad i \geq 0
\label{centralized-recursion}
}
where  $\mu>0$ is a small step-size parameter and for generality, we  are introducing  a diagonal matrix $D={\rm diag}\{d_{\ell} I_{M_\ell}\}_{\ell=1}^L \in \real^{M \times M}$, with $\{d_{\ell}\}$ being strictly positive (the choice $D=I$ is a special case). The initialization $w_{-1}$ is arbitrary. Since the gradient vector on the right-hand side of \eqref{centralized-recursion} is the sum of two separate gradient vectors shown in \eqref{centralized-gradient}, we can split the update into two incremental steps and write \cite{5,46}:
\begin{subequations}
\eq{
\psi_i &\ =\  w_{i-1} - \mu \eta D \ \grad_w p^{\rm glob}(w_{i-1}) \label{central-inc1}\\
w_i&\ =\  \psi_{i} - \mu D \ \grad_w J^{\rm glob}(\psi_i) \label{central-inc2}
}
\label{Centralized-sol}
\end{subequations}  
 \hspace{-2.5mm} In \eqref{central-inc1}--\eqref{central-inc2}, the vector $w_i$ is the estimate of the extended parameter $w$ at iteration $i$, while $\psi_i$ is an intermediate estimate for the same $w$. Note that the order of the incremental steps can be switched. However, in this work we consider penalizing first.
\subsection{Problem Reformulation for Distributed Solution}
In order to solve \eqref{penalized_cost} in a distributed manner, we first need to adjust the notation to account for one additional degree of freedom. Recall that the costs of two different agents, say, agents $k$ and $s$, may depend on the same sub-vector, say, $w^{\ell}$. Since these two agents will be learning $w^{\ell}$ over time, each one of them will have its own local estimate for $w^{\ell}$. While we expect these estimates to agree with each other over time, they nevertheless evolve separately and we need to use different notation to distinguish between them in the analysis. We do so by referring to the estimate of $w^{\ell}$ at agent $k$ by $w^{\ell}_k$ and to the estimate of $w^{\ell}$ at agent $s$ by $w^{\ell}_{s}$. In other words, we create virtual copies of the same sub-vector, $w^{\ell}$, with one copy residing at each agent. In this way, agent $k$ will evaluate iterates $w^{\ell}_{k,i}$ for $w^{\ell}$ over time $i$, and agent $s$ will evaluate iterates $w^{\ell}_{s,i}$ for the same $w^\ell$ over time $i$. As time evolves, we will show that these iterates will approach each other so that the subset of agents influenced by $w^\ell$ will reach agreement about it.
With this in mind, recall that we denoted the collection of all sub-vectors that influence agent $k$ by $w_k$; defined earlier in \eqref{w_I_k}. In that definition, the sub-vectors $\{w^{\ell}\}$ influencing $J_k(\cdot)$ were used  to construct $w_k$. In view of the new notation using virtual copies, we now redefine the same $w_k$ using the local copies instead, namely, we now write
\be
w_k \ \define\ \col\left\{w_k^{\ell}\right\}_{\ell \in{\cal I}_k}\;\in\real^{Q_k} \label{wk}
\ee
where $w_k^\ell \in \mathbb{R}^{M_\ell}$ is the local copy of the variable $w^\ell$ at agent $k$.  We  further let $\mathcal{C}_\ell$ denote the cluster of nodes that contains the variable $w^\ell$ in their costs:
\eq{\mathcal{C}_\ell= \{k \ |\ \ell \in \mathcal{I}_k\}
\label{cluster}}
We can view the cluster ${\cal C}_{\ell}$ as a smaller network (or sub-graph) where all agents in this sub-network are interested in the same parameter $w^{\ell}$. To require all local copies $w_{k}^{\ell}$ to coincide with each
other, we need to introduce the constraint
\be 
w_k^{\ell}
= w_{s}^{\ell},\;\;\forall \ k,s\in{\cal C}_{\ell}
\label{consensus}
\ee
Using relations \eqref{wk} and \eqref{consensus}, we rewrite problem \eqref{penalized_cost} as 
\begin{align}
 \underset{w_{1},....,w_{N}}{\text{minimize   }}& \quad
J^{\rm glob}_{\eta}(w_1,....,w_N) \triangleq \sum_{k=1}^N J_{k,\eta}(w_k) \label{penalized_cost2} \\
 \text{subject to     }& \quad
 w^\ell_k = w^\ell_s, \hspace{2 mm} \forall \hspace{2 mm} k, s \in \mathcal{C}_\ell, \hspace{1mm}\forall \hspace{1mm} \ell \nonumber
\end{align} 
\noindent The following example illustrates the above construction. \\

\noindent \begin{example}{\rm Consider the network with $5$ agents shown in Figure \ref{fig:illustrative(a)}.
\begin{figure}[H]
\centering
\begin{subfigure}[b]{0.5\textwidth}
\centering
\includegraphics[scale=0.55]{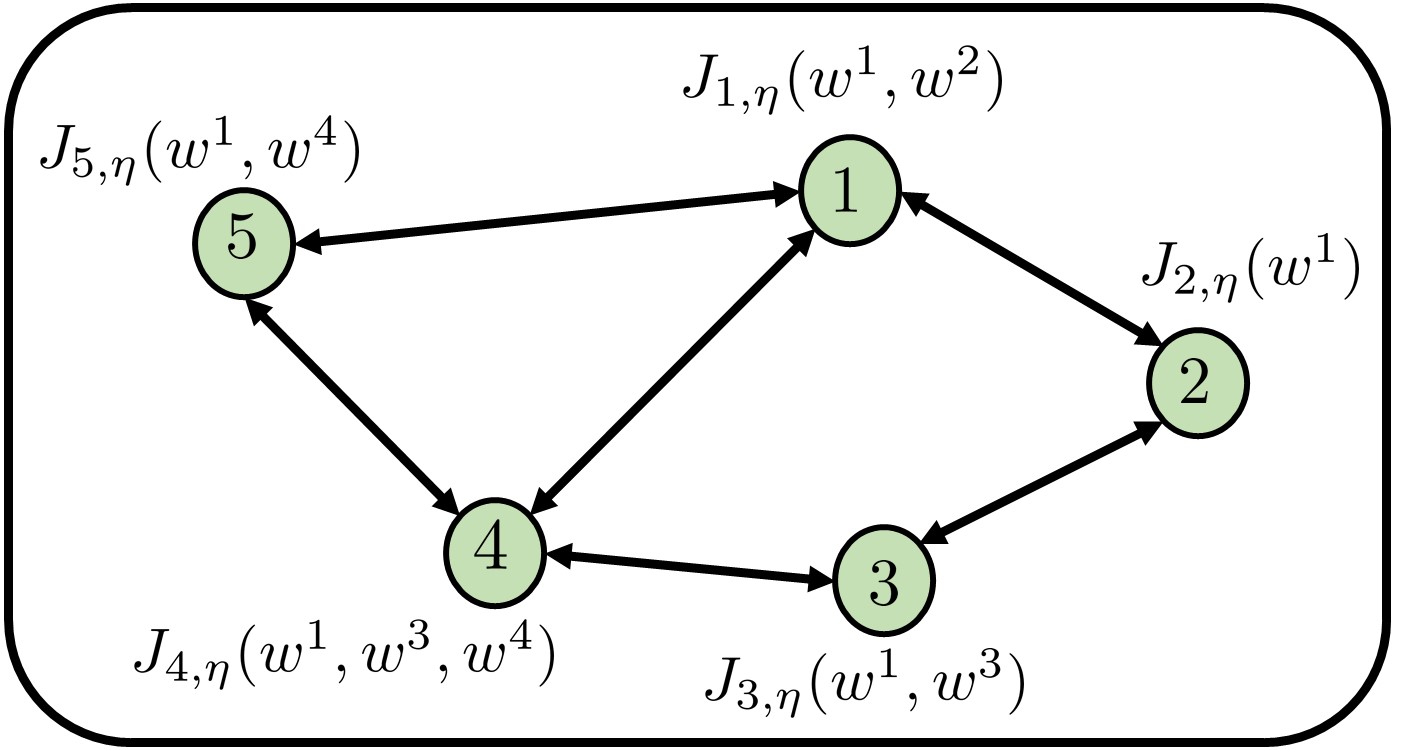}
\vspace*{-0.2cm}
\caption{}
\label{fig:illustrative(a)}
\end{subfigure}
\begin{subfigure}[b]{0.5\textwidth}
\centering
\includegraphics[scale=0.45]{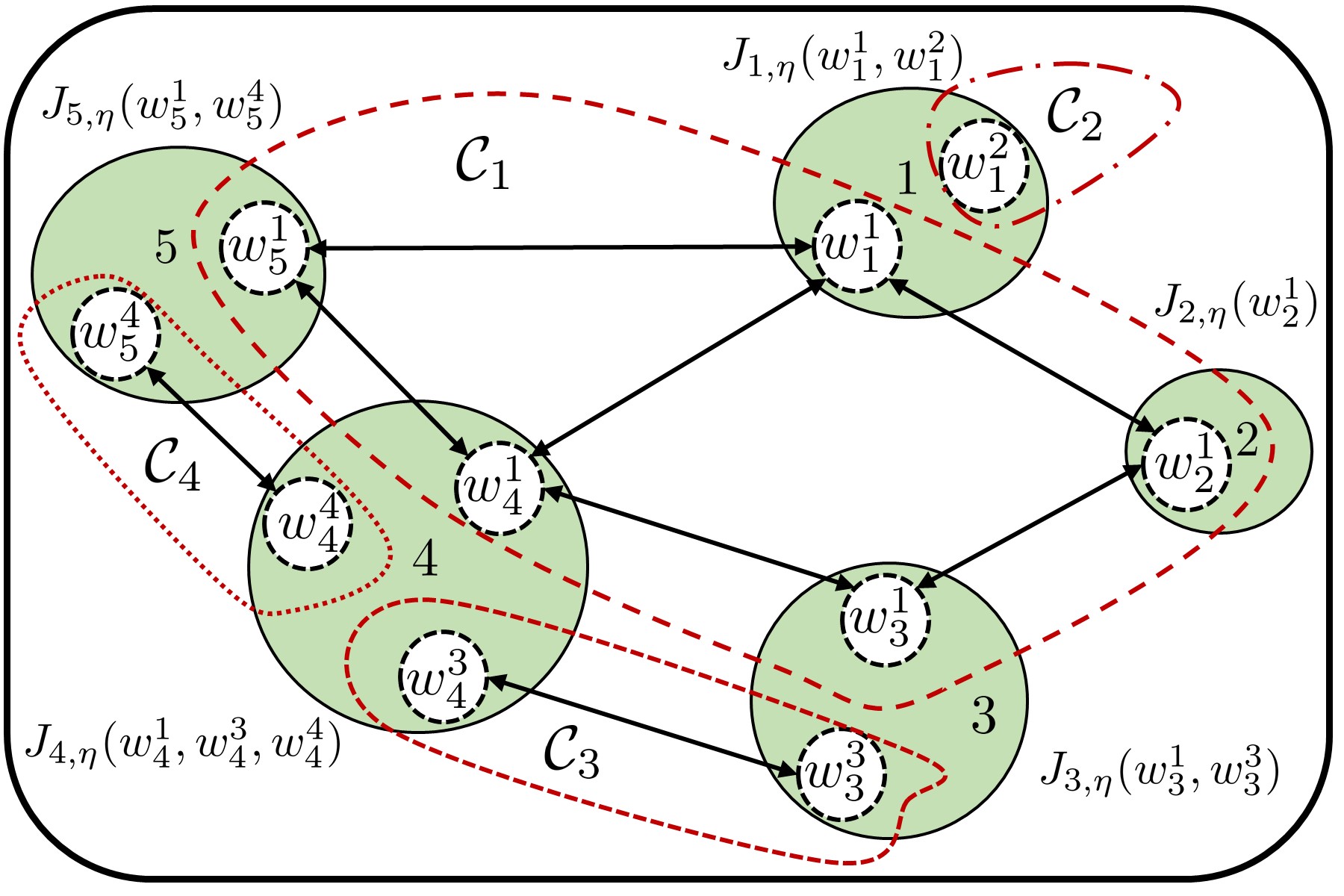}
\vspace*{-0.2cm}
\caption{}
\label{fig:illustrative(c)}
\end{subfigure}
\caption{A 5-agent network to illustrate the setting of problem \eqref{penalized_cost2}.}
\end{figure}

 In this network, we have $w = \col\{w^1, w^2, w^3,w^4\}$, $\cI_1 = \{1,2\}, \cI_2 = \{1\}$, $\cI_3=\{1,3\}$, $\cI_4=\{1,3,4\}$, and $\cI_5=\{1,4\}$. Consider further the penalized problem:
\eq{\label{xcwbsn}
\hspace{-2mm}\underset{\{w^1, w^2, w^3,w^4\}}{\mathrm{min.}}\, &J_{1,\eta}(w^1, w^2) + J_{2,\eta}(w^1) + J_{3,\eta}(w^1,w^3)+\nonumber \\
&J_{4,\eta}(w^1,w^3,w^4)+J_{5,\eta}(w^1,w^4)
}
To design a fully distributed algorithm, we introduce $w_k^{\ell}$ as the local copy of $w^{\ell}$ at agent $k$, and rewrite problem \eqref{xcwbsn} as:
\eq{\label{zw3nsdn}
{\text{minimize}}\quad &J_{1,\eta}(w_1^1, w_1^2) + J_{2,\eta}(w_2^1) + J_{3,\eta}(w_3^1,w_3^3)+\nonumber \\
&  J_{4,\eta}(w_4^1,w_4^3,w_4^4)+J_{5,\eta}(w_5^1,w_5^4), \nonumber\\
\text{subject\ to} \quad & w_1^1 = w_2^1=w_3^1=w_4^1=w_5^1 \nonumber \\
 \quad &w_3^3 = w_4^3 \nonumber \\
 \quad &w_4^4 = w_5^4 
}
If we introduce 
\eq{
w_1 &\define \col\{w_1^1, w_1^2\},\ w_2 \define \col\{w^1_2\},\ w_3 \define \col\{w_3^1, w_3^3\}, \nonumber \\
w_4 &\define \col\{w^1_4,w^3_4,w^4_4\},\ w_5 \define \col\{w_5^1, w_5^4\} 	\label{238sdh}
}
and organize the network into $L=4$ clusters as shown in Figure \ref{fig:illustrative(c)} with $\cC_1 = \{1,2,3,4,5\}$, $\cC_2=\{1\}$, $\cC_3 = \{3,4\}$, and $\cC_4 = \{4,5\}$, then problem \eqref{zw3nsdn} becomes equivalent to
\eq{
\underset{w_1, w_2, w_3,w_4}{\mbox{minimize}}& \quad \sum_{k=1}^{N} J_{k,\eta}(w_k), \nonumber \\
\mbox{subject to} &\quad w_k^\ell = w_s^\ell,\ \forall \hspace{1mm} k,s \in \cC_\ell,\ \ell = 1, 2, 3,4.
}
\qd}
\end{example}
\section{Coupled Diffusion Strategy}
\subsection{Cluster Combination Weights} \label{section_network_model}
 To solve \eqref{penalized_cost2}, we associate weights $\{a_{\ell,sk}\}_{s,k \in \cC_\ell}$ with each cluster $\cC_\ell$ and these weights are chosen to satisfy:
\eq{
&\sum_{s \in \cC_\ell} a_{\ell,sk} = 1, \quad a_{\ell,sk}=0, \ \text{if } s \notin \cN_k \cap \cC_\ell  \label{cluster_left_stochastic} 
}
\begin{assumption}
(\textrm{\bf{Each cluster is a connected sub-graph}}): \label{assump:connected}
The neighboring agents can communicate in both directions. Moreover, each ${\cal C}_{\ell}$ is connected. This implies that for any two arbitrary agents in cluster $\cC_\ell$, there exists at least one path with nonzero weights $\{a_{\ell,sk}\}_{s,k \in \cC_\ell}$ linking  one agent to the other. We also assume that at least one self weight $\{a_{\ell,kk}\}_{k \in \cC_\ell}$ is nonzero.			 			  \qd
	\end{assumption}
	\noindent 
We remark that two agents are coupled if they share the same variable $w^\ell$ and we are only requiring the coupled agents to be connected. If all agents share the entire $w$, then all agents are coupled by $w$ and the above assumption translates into requiring the network to be strongly-connected. Assumption \ref{assump:connected} is satisfied for most networks of interest. For example, all applications that fit into problem \eqref{glob-application} given in Remark \ref{remark-useful-case} naturally satisfy this assumption, including but not limited to applications in distributed power system monitoring, distributed control, and maximum-flow --- see  \cite{21,1,10}. This is due to the construction of problem \eqref{glob-application}: there exist $L=N$ clusters, where $w^\ell$ affects the neighborhood of agent $k=\ell$, and hence $\cC_k= \cN_k$ forms a star shaped graph (i.e., all agents $s \in \cC_k$ are connected through agent $k$), and hence this cluster is connected. Moreover, multitask applications satisfy this assumption \cite{4,14,15,16}. 

\begin{figure}[H]
	\centering
	\includegraphics[scale=.5]{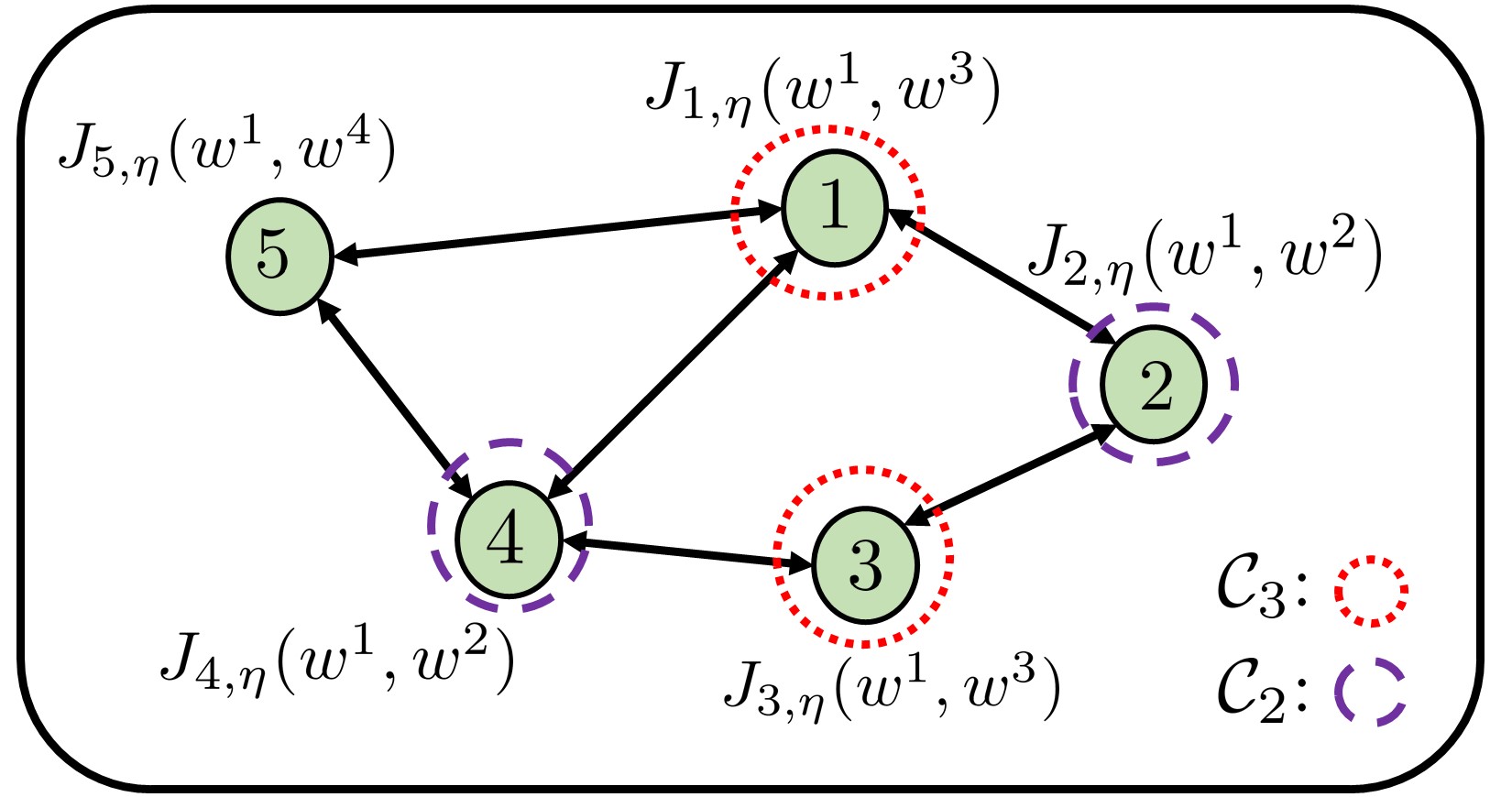}
	\caption{\small A five-agent network with unconnected $\cC_2$ and $\cC_3$.}
	\label{fig:illustrative-23hb}
\end{figure}

We emphasize that Assumption \ref{assump:connected} is not limited to the case described in Remark \ref{remark-useful-case} since it can 
be satisfied for any connected network, as we further clarify. To being with, independently of the clusters, let us assume that the entire network is connected. Now, if some cluster ${\cal C}_{\ell}$ happens to be unconnected, we can embed it into a larger {\em connected} cluster $\cC_{\ell}^\prime$ such that $\cC_{\ell} \subset \cC_{\ell}^\prime$. For example, consider the network shown in Figure \ref{fig:illustrative-23hb}.
In this network, we have
\eq{
\cC_1 = \{1,2,3,4,5\}, \ \cC_2 = \{2, 4\}, \ \cC_3 = \{1, 3\}, \ \cC_4=\{5\}
}
In these clusters, we find that $\cC_4$ is a singleton. Therefore, $w^4$ will be optimized solely and separately by  agent $5$, and no communication is needed for that variable. Cluster $\cC_1$ is connected, and agents $\{1,2,3,4,5\}$ cooperate in order to optimize $w^1$, with each agent sharing its estimate with neighbors. However, clusters ${\cal C}_2$ and $\cC_3$ have disconnected graphs. This implies that agents $2$ and $4$ cannot communicate to optimize and reach consensus on $w^2$. Likewise, for agents $\{1,3\}$ regarding the variable $w^3$. To circumvent this issue, we redefine $J_{1,\eta}(w^1, w^3)$ and $J_{2,\eta}(w^1,w^2)$ as:
\eq{
J_{1,\eta}'(w^1, w^2, w^3) &\define J_{1,\eta}(w^1, w^3) + 0 \cdot w^2 \\
J_{2,\eta}'(w^1, w^2, w^3) &\define J_{2,\eta}(w^1, w^2) + 0 \cdot w^3 
}
By doing so, the augmented  costs $J_{1,\eta}'(w^1, w^2, w^3)$ and $J_{2,\eta}'(w^1, w^2, w^3)$ now involve $w^2$ and $w^3$, respectively, and the new clusters become
\eq{
\cC_2^\prime = \{1, 2, 4\}, \quad \cC_3^\prime = \{1, 2, 3\}
}
which are connected and satisfy $\cC_2 \subset \cC_2^\prime$ and $\cC_3 \subset \cC_3^\prime$. Therefore, in this scenario, agents $\{1,2,4\}$ will now cooperate to optimize $w^2$ with agent $1$ acting as a connection that allows information about $w^2$ to diffuse in the cluster. Likewise, for agents $\{1,2,3\}$, with agent $2$ allowing information about $w^3$ to diffuse in the cluster. A second extreme approach would be to extend each local variable $w_k$ to the global variable $w$, which reduces problem \eqref{glob2} to the formulation \eqref{glob1}. This way of embedding the clusters into larger connected clusters can be done in a distributed fashion -- see for example \cite{12}.
\subsection{Coupled Diffusion Development}
Let $N_\ell$ denote the cardinality of cluster $\cC_\ell$ and further introduce the $N_\ell \times N_\ell$ matrix $A_{\ell}$ that collects the coefficients $\{a_{\ell,sk}\}_{s,k \in \mathcal{C}_{\ell}}$, namely,
\eq{ 
A_{\ell}\define [a_{\ell,sk}]_{s,k \in \mathcal{C}_{\ell}} \label{A_l}
} 
 Under condition \eqref{cluster_left_stochastic} and Assumption \ref{assump:connected}, the combination matrix $A_{\ell}$ will be left stochastic and primitive, i.e., $A_{\ell}\tran \one=\one$ and there exists a large enough $j_0$ such that the elements of $A_{\ell}^{j_0}$ are strictly positive. It follows from the Perron-Frobenius Theorem \cite[~ Lemma F.4]{9} that the matrix $A_{\ell}$ has a single eigenvalue at one of multiplicity one and all other eigenvalues are strictly less than one in magnitude. Moreover, the right eigenvector corresponding to the eigenvalue at one (the Perron vector), which we denote by $r_{\ell}$,
\eq{
A_{\ell}r_\ell=r_\ell, \quad  \one\tran r_\ell=1
}
is such that all its entries are positive and they are normalized to add up to one.

To facilitate the derivation that follows, we shall assume for the time being that $A_{\ell}$ is a locally balanced matrix \cite{45}, namely, that:
\eq{
R_\ell A_\ell\tran= A_\ell R_\ell \label{locally_balanced}
}
where $R_\ell \define {\rm diag}\{r_\ell(k)\}_{k \in \cC_\ell}$ is a diagonal matrix constructed from the Perron vector $r_\ell$ with $r_\ell(k)$ denoting the entry corresponding to agent $k \in \cC_\ell$. Condition \eqref{locally_balanced} is only used here to motivate the algorithm. Once derived, we will actually show that the algorithm is still convergent even when $A_\ell$ is only left stochastic but not necessarily locally balanced. To derive our proposed algorithm we state the following
auxiliary result proven in \cite{45}.
\begin{lemma} \label{null_v0}
	Let $R$ be a diagonal matrix constructed from the Perron vector of a left stochastic matrix $A \in \mathbb{R}^{Q \times Q}$. If $A$ is locally balanced matrix, i.e.,  $R A\tran= A R$, then it holds that $R- R A\tran$ is symmetric and positive semi-definite. Moreover, if we introduce the eigen-decomposition ${1\over2}(R- R A\tran)=U \Sigma U\tran$, the symmetric square-root matrix $Y\triangleq U \Sigma^{1/2} U\tran$ and let:
	\eq{
	\cR=R \otimes I_{M}, \quad	\sa=A \otimes I_{M}, \quad \mathcal{Y} = Y \otimes I_M
	}
	then, for primitive $A$ and any block vector $\ssx=\text{\em col}\{x^1,...,x^{Q}\}$ in the nullspace of $\cR-\cR \sa\tran $ with entries $x^q \in \mathbb{R}^{M}$ it holds that:
	\eq{
		\mathcal{Y}\ssx=0 \iff (\cR-\cR \sa\tran )\ssx=0  \iff x^1=x^2=...=x^{Q}		
	\label{nullspace}} \qd
\end{lemma} 
\noindent Lemma \ref{null_v0} allows us to rewrite \eqref{penalized_cost2} in an equivalent form that is amenable to distributed implementations. First, we introduce
\eq{\label{sadfads}
	\sw^{\ell} &\triangleq \text{col}\{w^{\ell}_k\}_{k \in \mathcal{C}_{\ell}} \in \real^{N_\ell M_\ell},
	}
which is the collection of all the local copies of $w^\ell$ across the agents in cluster $\mathcal{C}_\ell$. With this notation, we rewrite the cost function in problem \eqref{penalized_cost2} as
\eq{\label{equivalent-cost-function}
	\cJ(\sw^1, \sw^2, \cdots, \sw^L)+ \ \eta \cP(\sw^1, \sw^2, \cdots, \sw^L)
}
where
\eq{
\cJ(\sw^1, \sw^2, \cdots, \sw^L) &\define \sum_{k=1}^{N} J_{k}(w_k) \label{cost-cal-J} \\
\cP(\sw^1, \sw^2, \cdots, \sw^L) &\define \sum_{k=1}^{N} p_{k}(w_k) \label{cost-cal-P}
}
Next we use Lemma \ref{null_v0} to rewrite the constraints of problem \eqref{penalized_cost2} in an equivalent manner. We appeal to Lemma \ref{null_v0} to decompose
\eq{{1\over2} (R_\ell-R_\ell A_{\ell}\tran )=U_{\ell} \Sigma_{\ell} U_{\ell}\tran.
}
If we let
\eq{
	Y_{\ell} \triangleq U_{\ell} \Sigma_{\ell}^{1/2} U_{\ell}\tran, \quad \mathcal{Y}_{\ell} \triangleq  Y_{\ell} \otimes I_{M_{\ell}},
}
then using Lemma \ref{null_v0} and the definition of $\sw^\ell$ in \eqref{sadfads} we have 
\eq{\label{equi-constraint}
	w_k^\ell = w_s^\ell,\ \forall \ k,s \in \cC_\ell \Longleftrightarrow \cY_\ell \sw^\ell = 0,\quad \forall \hspace{1mm} \ell. 
}
Using relations \eqref{equivalent-cost-function} and \eqref{equi-constraint}, we can rewrite problem \eqref{penalized_cost2} equivalently as 
\begin{align}
\underset{\ssw^1,....,\ssw^L}{\text{minimize   }}& \quad  \cJ(\sw^1, \cdots, \sw^L)+\eta \cP(\sw^1, \cdots, \sw^L) \label{glob_exact0} \\
\text{subject to     }& \quad
\mathcal{Y}_{\ell}\sw^{\ell}=0, \hspace{1 mm} \forall \hspace{1 mm} \ell\nonumber
\end{align} 
To rewrite problem \eqref{glob_exact0} more compactly, we introduce
\eq{
\sw &\define \text{col}\{\sw^{\ell}\}_{\ell=1}^{L} \in \real^Q \label{wcaligraphic} \\
\mathcal{Y} &\define \text{blkdiag}\{\mathcal{Y}_{\ell}\}_{\ell=1}^L, \
	 \\
 \cJ(\sw) &\define \cJ(\sw^1, \cdots, \sw^L) \label{network-cost-J} \\
  \cP(\sw) &\define \cP(\sw^1, \cdots, \sw^L) \label{network-cost-P}  
}
where $Q\triangleq \sum\limits_{\ell=1}^L N_{\ell}M_{\ell}$. Then, problem \eqref{glob_exact0} becomes:
\begin{align}
\underset{\ssw}{\text{minimize }} \quad
\cJ(\sw)+\eta \cP(\sw),\hspace{1mm}  \label{glob_exact}
\text{s.t. } \mathcal{Y}\sw=0
\end{align} 
Instead of solving the constrained problem  \eqref{glob_exact}, we relax it and solve the penalized version:
\begin{align}
\boxed{
\underset{\ssw}{\text{minimize }} \quad
\cJ(\sw)+\eta \cP(\sw)+{1\over \mu}\|\cY \sw\|^2}  \label{glob_exact_penalized}
\end{align}  
with $\mu>0$. In \eqref{glob_exact_penalized} we see that ${1 \over \mu}$ is the penalty factor used for the consensus constraint \eqref{equi-constraint} and thus the smaller the value of $\mu$ is, the closer the solutions of problem \eqref{glob_exact} and \eqref{glob_exact_penalized} become to each other \cite{42,43}. We now note:
\eq{\cY^2=\text{blkdiag}\{\mathcal{Y}_{\ell}^2\}_{\ell=1}^L={1\over2}(\cR-\cR\sa\tran)
}
where
\eq{
\cR &\define {\rm blkdiag}\{\cR_{\ell}\}_{\ell=1}^L, \quad 
&\cR_\ell &\define R_\ell \otimes I_{M_\ell}
 \label{calblockR} \\
\sa &\define {\rm blkdiag}\{\sa_{\ell}\}_{\ell=1}^L, \quad 
&\sa_\ell &\define A_\ell \otimes I_{M_\ell}
 \label{calblockA}
}
  Applying three diagonally weighted incremental gradient descent steps to problem \eqref{glob_exact_penalized}, we get:
\begin{equation}
\begin{cases}
\begin{aligned}
 \zeta_i&=\sw_{i-1} - \mu \eta \cR^{-1} \grad_{\ssw} \cP(\sw_{i-1}) \\
 \psi_i&= \zeta_i -\mu \cR^{-1} \grad_{\ssw} \cJ(\zeta_{i}) \\
 \sw_i&=\psi_i - \mu \cR^{-1}\left({2\over \mu} \cY^2\right)\psi_i=\sa\tran \psi_i
\end{aligned}
\end{cases} 
\label{coupled-diff-network}
\end{equation} 
where $\mu$ is the step size.
Using the definition of $\cJ(\sw)$ and $\cP(\sw)$ from \eqref{network-cost-J}--\eqref{network-cost-P}, we have:
\eq{ &\scalemath{0.9}{
\grad_{\ssw} \cJ(\sw)= \begin{bmatrix}
\grad_{\ssw^1} \cJ(\sw) \\
\vdots \\
\grad_{\ssw^L} \cJ(\sw)
\end{bmatrix} },  \ \scalemath{0.9}{ \grad_{\ssw^\ell} \cJ(\sw)} = \text{col}\{\grad_{w_k^\ell} J_k(w_k)\}_{k \in \mathcal{C}_\ell} \label{gradient-network-cost}\\
&\scalemath{0.9}{ \grad_{\ssw} \cP(\sw)= \begin{bmatrix}
\grad_{\ssw^1} \cP(\sw) \\
\vdots \\
\grad_{\ssw^L} \cP(\sw)
\end{bmatrix} }, \ \scalemath{0.9}{ \grad_{\ssw^\ell} \cP(\sw)} = \text{col}\{\grad_{w_k^\ell} p_k(w_k)\}_{k \in \mathcal{C}_\ell} \label{gradient-network-penalty}
}
Therefore, using the definition of $\sa$ from \eqref{calblockA} and $R_\ell = {\rm diag}\{r_\ell(k)\}_{k \in \cC_\ell}$, recursion \eqref{coupled-diff-network} can be rewritten more explicitly in distributed form as listed below.
\begin{subequations}
\label{Coupled diffusion}
\eq{
\zeta_{k,i}&=w_{k,i-1}- \mu \eta \Omega_k \grad_{w_k} p_k(w_{k,i-1}) \label{Coupled diffusion-(a)}\\
\psi_{k,i}&=\zeta_{k,i}-\mu \Omega_k \grad_{w_k} J_k(\zeta_{k,i}) \label{Coupled diffusion-(b)}\\
 w^\ell_{k,i}&= \sum_{\substack{s \in \mathcal{N}_k \cap \mathcal{C}_\ell}} a_{\ell,sk}\psi^\ell_{s,i}, \hspace{2mm} \forall \hspace{1mm} \ell \in \cI_k \label{Coupled diffusion-(c)}
}
\end{subequations}
where 
\eq{
\Omega_k={\rm diag} \left\{ I_{M_\ell}/ r_\ell(k) \right\}_{\ell \in \cI_k} \label{Omega_step-size}
}
and $r_\ell(k)$ denote the entry of the Perron vector $r_\ell$ corresponding to agent $k \in \cC_\ell$. In this description, the variables $\{\zeta_{k,i}, \psi_{k,i}\}$ are intermediate estimates for the vector $w_k$, which contains all the parameters that influence agent $k$. These vectors have dimension $Q_k\times 1$ each. On the other hand, the variable  $w_{k,i}^{\ell}$ has size $M_{\ell}\times 1$  and is an estimate for the specific parameter of index $\ell$, i.e., $w_k^{\ell}$. Thus, note that in  steps \eqref{Coupled diffusion-(a)}--\eqref{Coupled diffusion-(b)}, a traditional diagonally weighted gradient-descent step is applied by each agent using the gradients of the corresponding penalty and risk functions; these steps generate the intermediate $Q_k-$dimensional iterates $\{\zeta_{k,i},\psi_{k,i}\}$. The last step \eqref{Coupled diffusion-(c)} is a convex combination step, where each agent $k$ combines the iterates of index $\ell$ from their neighbors to construct $w_{k,i}^{\ell}$. More specifically, for every $\ell \in \cI_k$, each agent $k$ combines its entry $\psi_{k,i}^{\ell}$ with the neighboring entries $\{\psi_{s,i}^\ell \ | \ s \in \cN_k \cap \cC_\ell\}$ using weights $\{a_{\ell,sk}\}_{s \in \cN_k \cap \cC_\ell}$.
It should be noted that each agent $k$ gets to choose its own combination weights. For example, let $n_{\ell,k}= \big|\mathcal{N}_k \cap \mathcal{C}_\ell \big|$ denote the number of agents that belong to $\cC_\ell$ and are neighbors of agent $k$ (including agent $k$). Then, we can use the Metropolis rule to construct the combinations weights $\{a_{\ell,sk}; \ s \in \cN_k \cap \mathcal{C}_\ell,\ \ell \in \cI_k\}$ as follows \cite{9}:
\begin{equation}
a_{\ell,sk} =
\begin{cases}
\begin{aligned}
  &{1 \over \max\{n_{\ell,k},n_{\ell,s}\}}
,& \quad& \text{if }  s \in \mathcal{N}_k \cap \mathcal{C}_\ell,\ s \neq k  \\
&1-\sum_{e \in \cN_k \cap \mathcal{C}_\ell \backslash \{k\}} a_{\ell,ek},& \quad& s=k,  \\
&0, & \quad& \text{otherwise.} 
\end{aligned}
\end{cases} \label{Metropolis}
\end{equation}
or we can use the averaging rule \cite{9}:
\begin{equation}
a_{\ell,sk} =
\begin{cases}
\begin{aligned}
  &{1 \over n_{\ell,k}}
,& \quad& \text{if }  s \in \mathcal{N}_k \cap \mathcal{C}_\ell,  \\
&0, & \quad& \text{otherwise.} 
\end{aligned}
\end{cases} \label{averaging_rule}
\end{equation}
 \noindent \begin{remark} {\rm \textbf{(Step-size)}\label{remark-step-size} Due to the use of left stochastic matrices, and in order for the algorithm to converge to the minimizer of the sum of costs, the step-sizes used for each entry $w^\ell_{k,i-1}$ in $w_{k,i-1}$ need to be divided by the entry of the Perron eigenvector  $r_\ell(k)$ corresponding to cluster $\ell$ for agent $k$ as in \eqref{Coupled diffusion-(a)}--\eqref{Coupled diffusion-(b)}. Otherwise, the algorithm will converge to a neighborhood of a different limit point that satisfies:
\eq{
\sum_{k \in \mathcal{C}_\ell} r_\ell(k) \grad_{w_k^\ell}J_{k,\eta}(w_k)=0 , \quad \forall \ \ell
} 
For the case of $L=1$ and $\cI_k=\{L\}$, this will correspond to a Pareto optimal solution \cite{7,9}. We further remark that for many combination rules, the entries $r_\ell(k)$ are known. For example, for the Metropolis rule $r_\ell(k)=(N_\ell)^{-1}$, while for the averaging rule it is $r_\ell(k)=n_{\ell,k} \left(\sum_{s \in \cC_\ell}n_{\ell,s}\right)^{-1}$. This is not an issue since the common part can be absorbed into $\mu$. For example, for the averaging rule we have $\mu /r_\ell(k)  =\mu' n_{\ell,k}^{-1} $, where now the agents need to agree on $\mu'=\mu  \left(\sum_{s \in \cC_\ell}n_{\ell,s}\right)$ instead of $\mu$. For general left stochastic matrices, results already exist in the literature that can estimate the Perron entries in a distributed fashion \cite{45}.
\qd }
\end{remark}
\section{Stochastic Analysis Setup}  
As mentioned before in Remark \ref{remark-special-case}, we will also allow for the possibility of stochastic risks, in which case the true gradient vectors are not available. Therefore, we introduce the gradient noise vector for each agent at time $i$:
\eq{
  \v_{k,i}(\bzeta_{k,i})  & \define   \grad_{w_k} J_k(\bzeta_{k,i}) -\grad_{w_k} Q(\bzeta_{k,i};\x_{k,i}) \nonumber \\
&  \hspace{1mm}= \grad_{w_k} J_k( \bzeta_{k,i}) - \widehat{\grad_{w_k} J}_k( \bzeta_{k,i})  \label{gradient-model} 
 }
  that is required to satisfy certain conditions given in Assumption \ref{noisemodel:assump}.
\begin{assumption} \label{noisemodel:assump}
(\textrm{\bf{Gradient noise model}}): Conditioned on the past history of iterates $\cf_{i} \triangleq \{ \w_{k,j-1} : k=1 , ... , N \text{ and } j \leq i \}$, the gradient noise $\v_{k,i}(\bzeta_k)$ is assumed to satisfy:
\eq{
\Ex\{\v_{k,i}(\bzeta_k)\mid \cf_{i} \} &= 0 \label{noise-model(a)} \\
\Ex\{\|\v_{k,i}(\bzeta_k)\|^2\mid \cf_{i} \} &\leq  \bar{\alpha}_k \|\bzeta_k\|^2 + \bar{\sigma}_k^2 \label{noise-model(b)}
} 
for some nonnegative constants $\bar{\alpha}_k$ and $\bar{\sigma}_k^2$. 
 \qd
\end{assumption}
We again emphasize that in this work we account for noisy gradients like \eqref{gradient-model} and, therefore, we shall incorporate the presence of the gradient noise into the analysis. In the presence of stochastic gradient constructions, the
coupled diffusion algorithm \eqref{Coupled diffusion} becomes the one listed in \eqref{GM Stochastic Diffusion}. Note that we are now using boldface letters in \eqref{GM Stochastic Diffusion} to highlight the fact that the variables are stochastic in nature due to the randomness in the gradient noise component.
\begin{algorithm}[H]
\caption{(Coupled diffusion strategy)}
{\bf Setting:} Let $\Omega_k={\rm diag} \left\{I_{M_\ell} / r_\ell(k) \right\}_{\ell \in \cI_k}$ and $\w_{k,-1}$ arbitrary. \\
{\bf For every agent $k$, repeat for $i\geq 0$:}
\begin{subequations}
\label{GM Stochastic Diffusion}
\eq{
\bzeta_{k,i}&= \w_{k,i-1}-\mu\eta \Omega_k  \grad_{w_k}  p_k(\w_{k,i-1}) \label{stochastic-diff(a)}\\
\bpsi_{k,i}&=\bzeta_{k,i}-\mu \Omega_k  \widehat{\grad_{w_k} J}_k( \bzeta_{k,i}) \label{stochastic-diff(b)}\\ 
 \textrm{\bf For}& \ \textrm{\bf every block entry  $\ell \in \cI_k$, combine:}  \nonumber \\
& \w^\ell_{k,i}= \sum_{\substack{s \in \mathcal{N}_k \cap \mathcal{C}_\ell}} a_{\ell,sk}\bpsi^\ell_{s,i} \label{stochastic-diff(c)} 
}
\end{subequations}
\end{algorithm}

We will measure the performance of the distributed strategy by examining the mean-square-error between the random iterates $\boldsymbol \w^\ell_{k,i}$ and the corresponding optimal component from \eqref{optimal-original}, denoted by $w^{\ell,o}$. For this purpose, we first note that in terms of the optimal solution $w^{\ell,\star}$ for the penalized problem \eqref{optimal-penalized}, we can write:
\eq{
&\limsup  \limits_{i\rightarrow \infty}\hspace{1mm} \Ex \|w^{\ell,o} - \w^\ell_{k,i} \|^2 \nonumber \\
& \quad =\limsup \limits_{i\rightarrow \infty}
\Ex \|w^{\ell,o} -w^{\ell,\star}+w^{\ell,\star}- \w^\ell_{k,i} \|^2 \nonumber \\
&\quad \leq 2 \underbrace{\|w^{\ell,o} -w^{\ell,\star}\|^2}_{\text{Approximation Error}} +2 \limsup \limits_{i\rightarrow \infty} \Ex \|w^{\ell,\star}- \w^\ell_{k,i} \|^2 \label{difference to optimal}
} 
The following result is proven in \cite{5} concerning the size of the first component as the size of the penalty factor becomes unbounded.
\begin{theorem}(\textrm{\bf{Approaching Optimal Solution}}): Under Assumptions \ref{feasible assump}--\ref{penalty-assump}, it holds that:
\eq{
\lim  \limits_{\eta \rightarrow \infty}
\|w^o -w^\star\| = 0
\label{th-approachoptiml}} \label{theorem-approachoptimal}
\qd
\end{theorem}
\noindent Therefore, to assess \eqref{difference to optimal} we will characterize the second term $\Ex \|w^{\ell,\star}- \w^\ell_{k,i} \|^2$ and show that we can drive it to arbitrarily small values.
In order to carry out the analysis we need to examine the error dynamics of the algorithm more closely. 

For ease of reference we collect all the main symbols into the following table.
\begin{table}[h]
	\centering \caption{A listing of the main symbols and their interpretation.}
	\begin{tabular}{ | c || l | }
		\hline  \hline
		\cellcolor{gray!25} \bf Symbol &   \multicolumn{1}{c|}{ \bf \cellcolor{gray!25} Meaning} \\
	 \hline
		 $\cI_k$ & The set of variable indices that influence the cost of agent $k$. 
		\\  \hline
		$w^\ell_{k}$ & Local copy of $w^\ell$ at agent $k$. \\
		 \hline
		$w_{k}$ & Stacks the parameters influencing agent $k$,  $w_{k} \triangleq  \text{col}\{w^{\ell}_k\}_{\ell \in \mathcal{I}_k}$ \\
		 \hline
		$\mathcal{C}_\ell$ & Cluster of nodes that is influenced by the variable $w^\ell$.\\ 
		\hline 
		$\sw^\ell$ & Stacks all local copies of $w^{\ell}$ across ${\cal C}_{\ell}$,
 $\sw^{\ell} \triangleq \text{col}\{w^{\ell}_k\}_{k \in \mathcal{C}_{\ell}}$  \\
		\hline
		$\sw$ & Stacks $\sw^{\ell}$ for all parameters, $\sw\triangleq \text{col}\{\sw^{\ell}\}_{\ell=1}^L$\\
		\hline
		$\cJ(\sw)$ & Global risk, $\cJ(\sw)\triangleq \sum_{k=1}^N J_k(w_k)$  \\
		\hline
	$\cP(\sw)$ & Global penalty, $\cP(\sw)\triangleq \sum_{k=1}^N p_k(w_k)$\\  
		\hline \hline
	\end{tabular}
	\label{table-notation}
\end{table}
\section{Error Dynamics}
\subsection{Network Error Recursion}
We start by expanding \eqref{gradient-model} into its individual components:
\eq{
 \v_{k,i}^\ell(\bzeta_{k,i})   =   \grad_{w_k^\ell} J_k( \bzeta_{k,i}) -\widehat{\grad_{w_k^\ell} J}_k( \bzeta_{k,i}) , \ \ell \in \cI_k \label{gradient-model-components}
 }
where $\v_{k,i}^\ell(\bzeta_{k,i})$ is the part of the gradient noise related to approximating $\grad_{w_k^\ell} J_k(\bzeta_{k,i})$. We collect the noise terms $\{\v_{k,i}^\ell(\bzeta_{k,i})\}_{k \in \cC_\ell}$ across all agents and clusters into block vectors:
\eq{
 \v_{i}^\ell \define {\rm col}\{\v_{k,i}^\ell(\bzeta_{k,i})\}_{k \in \cC_\ell}, \quad 
\v_{i} \define {\rm col}\left\{\v_i^\ell \right\}_{\ell=1}^L  \label{transform-vector-noise}
}
Similarly, motivated by \eqref{sadfads} and \eqref{wcaligraphic} we define the network vectors:
\eq{
\bsw_i^{\ell} \define \text{col}\{\w^{\ell}_{k,i}\}_{k \in \mathcal{C}_{\ell}}, \quad \bsw_i \define \text{col}\{\bsw_i^{\ell}\}_{\ell=1}^{L},
}
 Incorporating the gradient noises \eqref{transform-vector-noise} into \eqref{coupled-diff-network} we obtain the network recursion of \eqref{GM Stochastic Diffusion}:
 \begin{subequations}
\eq{
\bzeta_i&=\bsw_{i-1} - \mu \eta \cR^{-1}\grad_{\ssw} \cP(\bsw_{i-1}) \label{network-recur-noise(a)}\\
 \bpsi_i&= \bzeta_i -\mu \cR^{-1} \grad_{\ssw} \cJ(\bzeta_i)+\mu \cR^{-1} \v_i \label{network-recur-noise(b)}\\
 \bsw_i&=\sa\tran \bpsi_i \label{network-recur-noise(c)}
 }
  \end{subequations}
Now recall that problems \eqref{penalized_cost} and \eqref{glob_exact} are equivalent and, therefore, the optimal solution to \eqref{glob_exact} is given by
\eq{
\sw^\star \define {\rm col}\{\sw^{\ell,\star}\}_{\ell=1}^L, \quad
\sw^{\ell,\star} \define \one_{N_\ell} \otimes w^{\ell,\star}
}
Subtracting $\sw^\star$ from both sides of \eqref{network-recur-noise(a)}--\eqref{network-recur-noise(c)} we get:
\begin{subequations}
\eq{
\tzeta_i&=\tsw_{i-1} + \mu \eta \cR^{-1} \grad_{\ssw} \cP(\bsw_{i-1}) \label{zeta-error}
\\
\tpsi_{i} &= \tzeta_{i} + \mu \cR^{-1} \grad_{\ssw} \cJ(\bzeta_i) - \mu \cR^{-1} \v_i \label{psi-error} \\
\tsw_{i} &= \sa\tran \tpsi_{i} \label{Werror}
}
\end{subequations}
where $\tsw_{i-1}\define \sw^\star-\bsw_{i-1}$ denotes the error at time $i-1$ and similarly for $\{\tzeta_i,\tpsi_{i}\}$.
Combining \eqref{zeta-error}, \eqref{psi-error}, and \eqref{Werror} we arrive at the following statement.
\begin{lemma}(\textrm{\bf{Network error recursion}}). The network error vector evolve according to the following dynamics:
\eq{
\scalemath{0.95}{ \tsw_{i}= \sa\tran \bigg(\tsw_{i-1} + \cR^{-1} \big( \mu \eta  \grad_{\ssw} \cP(\bsw_{i-1})+ \mu  \grad_{\ssw} \cJ(\bzeta_i)  - \mu  \v_i \big) \bigg) } \label{error1}
}
where
\eq{
\scalemath{0.95}{ \bzeta_i=\bsw_{i-1} - \mu \eta \cR^{-1} \grad_{\ssw} \cP(\bsw_{i-1})  }
}
\qd
\end{lemma}
\section{Transformed Network Error Dynamics}
\subsection{Similarity Transformation I}
The convergence analysis of recursion \eqref{error1} is facilitated by transforming it to a convenient basis. Since each $A_{\ell}$ is left-stochastic and primitive, it admits a Jordan decomposition of the form \cite{9}:
\eq{
A_\ell &\triangleq  V_{\ell} J_\ell V_{\ell}^{-1}
}
where
\eq{
  V_{\ell} = \begin{bmatrix}
 r_\ell  & V_{1,\ell} \end{bmatrix}, \ \  J_\ell  &= \begin{bmatrix}
 1 & 0 \\
 0 & \check{J}_{\ell}
 \end{bmatrix} , \ \ 
 V_{\ell}^{-1} = \begin{bmatrix}
 \one_{N_\ell}\tran \\
 V_{2,\ell}\tran
 \end{bmatrix} \label{decomposition-A_l} 
}
and the matrices $V_{1,\ell}$ and $V_{2,\ell}$ have dimensions $N_\ell \times (N_\ell-1)$. Moreover, the matrix $\check{J}_{\ell}$ is $(N_\ell-1) \times (N_\ell-1)$ and consists of Jordan blocks, with each one having the generic form (with a positive scalar $\epsilon$ replacing the usual unit entries in the lower diagonal):
 \eq{ \scalemath{0.9}{
\begin{bmatrix}
\lambda & &\\
\epsilon & \lambda  & \\
&   \ddots & \ddots & \\
&  & \epsilon& \lambda
\end{bmatrix}}
}
If we define
\eq{ \mathcal{V}_{\ell} &\triangleq V_{\ell} \otimes I_{M_\ell}=\begin{bmatrix}
r_\ell \otimes I_{M_\ell} & V_{1,\ell}\otimes I_{M_\ell} 
\end{bmatrix} \label{calV_ell}\\
 \mathcal{J}_\ell &\triangleq  J_\ell \otimes I_{M_\ell} = \begin{bmatrix}
 I_{M_\ell} & 0 \\
 0 & \check{J}_{\ell} \otimes I_{M_\ell}
 \end{bmatrix} \label{calJ_ell}
}
then we can decompose $\sa$ from \eqref{calblockA} as:
\eq{
\sa =&  \underbrace{\text{ blkdiag}\{\mathcal{V}_{\ell} \}_{\ell=1}^L}_{\define \cV} \underbrace{\text{ blkdiag}\{\mathcal{J}_\ell \}_{\ell=1}^L}_{\define \mathcal{J}} \underbrace{\text{ blkdiag}\{\mathcal{V}_{\ell}^{-1} \}_{\ell=1}^L}_{\define \cV^{-1}}
}
We now multiply both sides of the error recursion \eqref{error1} from the left by $\cV\tran$:
\eq{
\cV\tran \tsw_{i}&=    \cJ\tran \bigg( \cV\tran \tsw_{i-1} + \mu \eta \cV\tran \cR^{-1} \grad_{\ssw} \cP(\bsw_{i-1}) \nonumber \\
& \quad + \mu \cV\tran \cR^{-1} \grad_{\ssw} \cJ(\bzeta_i) - \mu \cV\tran \cR^{-1} \v_i\bigg)  
 \label{trans-error} 
 }
and denote the individual block entries of the various quantities in \eqref{trans-error} by
\begin{align}
\cV\tran \tsw_{i} ={\rm col}\left\{\begin{bmatrix}
(r_\ell\tran \otimes I_{M_\ell}) \tsw_{i}^\ell \\ 
(V_{1,\ell}\tran \otimes I_{M_\ell}) \tsw_{i}^\ell
 \end{bmatrix}\right\}_{\ell=1}^L 
 &\triangleq  {\rm col}\left\{\begin{bmatrix}
\bar{\w}_{i}^\ell \\ 
\cw_{i}^\ell
 \end{bmatrix} \right\}_{\ell=1}^L \label{w-cal-bar-check} \\
 \mu \cV\tran \cR^{-1} \grad_{\ssw} \cJ(\bzeta_i)
 &\triangleq  {\rm col}\left\{\begin{bmatrix}
\bar{\g}^\ell(\bzeta_i) \\ 
\check{\g}^{\ell}(\bzeta_i)
 \end{bmatrix} \right\}_{\ell=1}^L \label{g-bar-check} \\
   \mu \eta \cV\tran \cR^{-1} \grad_{\ssw} \cP(\bsw_{i})  & \triangleq  {\rm col}\left\{\begin{bmatrix}
\bar{\f}^\ell(\bsw_{i}) \\ 
\check{\f}^\ell(\bsw_{i})
 \end{bmatrix} \right\}_{\ell=1}^L   \label{f-bar-check}
 \\
 -\mu  \cV\tran   \cR^{-1}  \v_i &\triangleq  {\rm col}\left\{\begin{bmatrix}
\bar{\v}_{i}^\ell \\ 
\check{\v}_{i}^\ell
 \end{bmatrix} \right\}_{\ell=1}^L \label{v-bar-check} 
\end{align}
 Note that the quantities $\{\bar{\w}_{i}^\ell,\bar{\g}^\ell(\bzeta_{i}),\bar{\f}^\ell(\bsw_{i}),\bar{\v}_{i}^\ell\}$ are vectors of dimension $M_\ell \times 1$ while the quantities $\{\cw_{i}^\ell,\check{\g}^\ell(\bzeta_{i}),\check{\f}^{\ell}(\bsw_{i}),\check{\v}_{i}^\ell\}$ are vectors of dimension  $M_\ell (N_\ell-1) \times 1$. Using these transformations, we can rewrite the previous recursion \eqref{trans-error} as:
 \eq{
\scalemath{0.96} { {\rm col}\left\{\begin{bmatrix}
\bar{\w}_{i}^\ell \\ 
\cw_{i}^\ell
 \end{bmatrix} \right\}  =  \mathcal{J}\tran \ {\rm col}\left\{\begin{bmatrix}
\bar{\w}_{i-1}^\ell+\bar{\g}^\ell(\bzeta_{i})+\bar{\f}^\ell(\bsw_{i-1})+\bar{\v}_{i}^\ell \\ 
\cw_{i-1}^\ell+\check{\g}^\ell(\bzeta_{i})+\check{\f}^{\ell}(\bsw_{i-1})+\check{\v}_{i}^\ell
 \end{bmatrix} \right\} }
\label{errorJ}
}
With a slight abuse of notation we dropped $\ell$ from ${\rm col}\{.\}_{\ell=1}^L$ and wrote it as ${\rm col}\{.\}$, we continue to do so in the following unless it is not clear from the presentation.
\subsection{Similarity Transformation II}
\begin{figure}[H]
\centering
\includegraphics[scale=0.3]{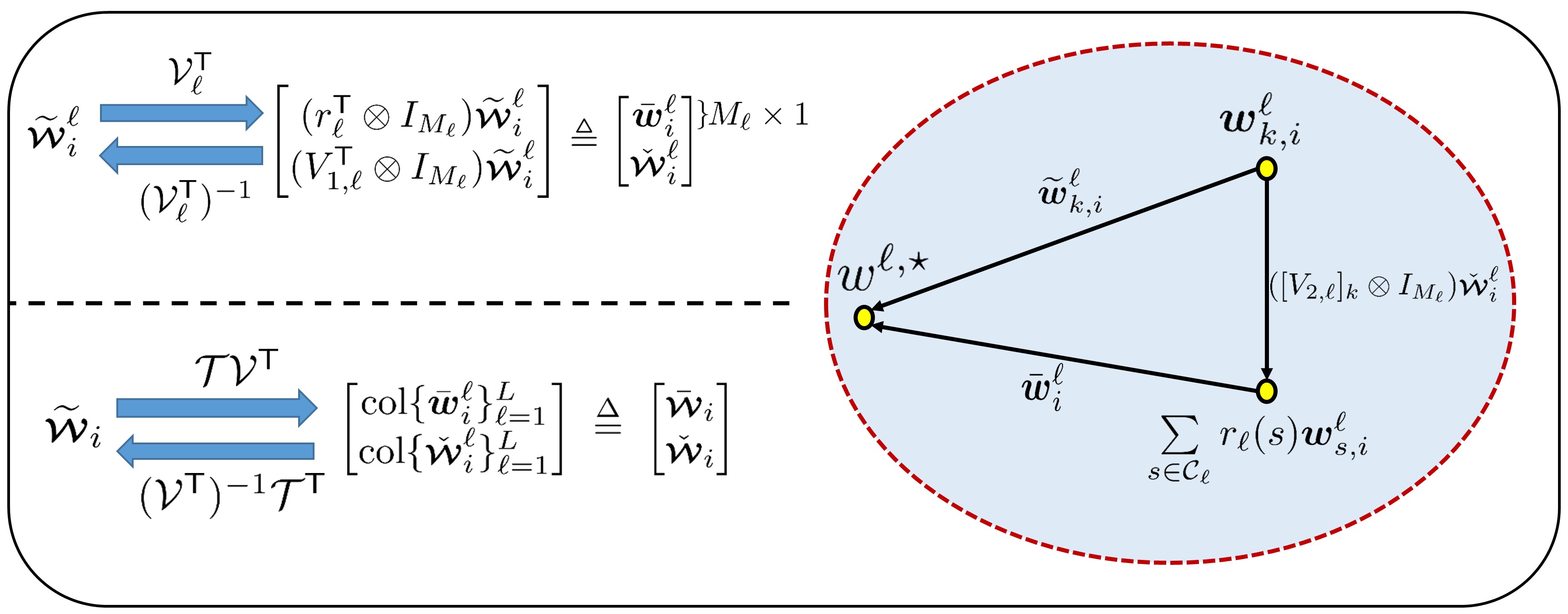}
\caption{Geometric relation between the error $\tw^\ell_{k,i}$ and the transformed parts $\bar{\w}_{i}^\ell$ and $\cw_{i}^\ell$.}
\label{fig:relationerr}
\end{figure}

\begin{figure}[H]
\includegraphics[scale=0.28]{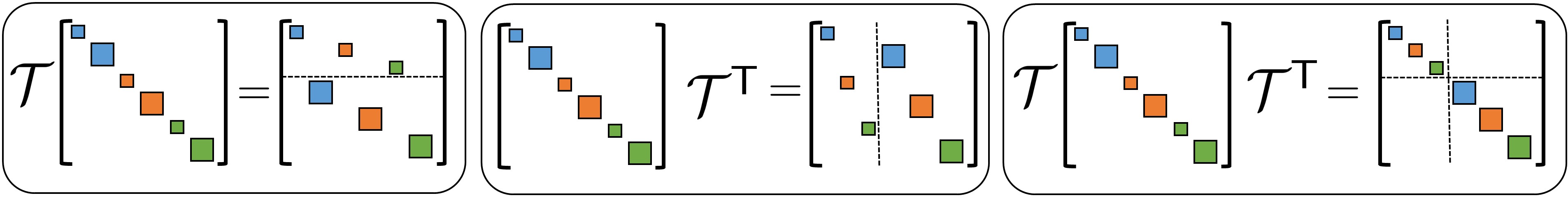}
\caption{A visual illustration of transformations $\st \mathcal{J}\tran \st\tran$ for $L=3$ blocks represented by different colors.}
\label{fig:TXT}
\end{figure}

If we examine the structure of the entries in the transformed vector \eqref{w-cal-bar-check}, we observe that 
\eq{
\bar{\w}_{i}^\ell =\sum_{k \in \cC_\ell} r_\ell(k)\tw_{k,i}^\ell \label{centroid_error}
}
 is the centroid of the errors relative to $w^{\ell,\star}$ across all agents in $\cC_\ell$. Moreover, we note that
\eq{
\tsw_i^\ell=(\cV_\ell^{-1})\tran \begin{bmatrix}
\bar{\w}_{i}^\ell \\ 
\cw_{i}^\ell
 \end{bmatrix}  =  \one_{N_\ell} \otimes \bar{\w}_{i}^\ell +  (V_{2,\ell} \otimes I_{M_\ell}) \cw_{i}^\ell \label{clustererr=avg+dev}
}
so that we can write each $\tw_{k,i}^\ell$ as
\eq{
\tw_{k,i}^\ell  =   \bar{\w}_{i}^\ell +  \left([V_{2,\ell}]_k \otimes I_{M_\ell} \right) \cw_{i}^\ell
\label{erro-ind-av+dev}}
where $\left([V_{2,\ell}]_k \otimes I_{M_\ell} \right)$ are the rows of $\left(V_{2,\ell} \otimes I_{M_\ell} \right)$ corresponding to the position of $\tw_{k,i}^\ell$ in the cluster vector $\tsw_i^\ell$.  From \eqref{erro-ind-av+dev} we see that $\{([V_{2,\ell}]_k \otimes I_{M_\ell}) \cw_{i}^\ell\}$ is the deviation of the individual errors $\{\tw_{k,i}^\ell\}$ from the weighted centroid error $\bar{\w}_{i}^\ell$ -- see Figure \ref{fig:relationerr}. 
In the following we will show that we can drive $\bar{\w}_{i}^\ell$ and $\cw_{i}^\ell$ to arbitrarily small values. We first re-order the elements in \eqref{w-cal-bar-check}, so that iterates that correspond to the weighted centroids $\{\bar{\w}_{i}^\ell\}$ appear stacked together. For this purpose, we introduce a permutation matrix ${\cal T}$ such that
\eq{\st {\rm col}\left\{\begin{bmatrix}
\bar{\w}_{i}^\ell \\ 
\cw_{i}^\ell
 \end{bmatrix} \right\}_{\ell=1}^L =\begin{bmatrix}
{\rm col}\{
\bar{\w}_{i}^\ell
  \}_{\ell=1}^L \\ 
{\rm col}\{ 
\cw_{i}^\ell \}_{\ell=1}^L
 \end{bmatrix} \define \begin{bmatrix}
\bw_i \\ 
\cw_i
 \end{bmatrix} \label{permutaion}
 }
 We then transform the error recursion \eqref{errorJ} by multiplying both sides by $\mathcal{T}$ on the left to get:
 \eq{
\begin{bmatrix}
\bw_i \\ 
\cw_i
 \end{bmatrix}= & \st \mathcal{J}\tran \st\tran  \begin{bmatrix}
\bw_{i-1}+\bar{\g}(\bzeta_{i})+\bar{\f}(\bsw_{i-1})+\bar{\v}_{i} \\ 
\cw_{i-1} +\check{\g}(\bzeta_{i})+\check{\f}(\bsw_{i-1})+\check{\v}_{i}
 \end{bmatrix}  
\label{errorT}} 
where 
\eq{
\begin{bmatrix}
\bar{\g}(\bzeta_{i}) \\ 
\check{\g}(\bzeta_{i})
 \end{bmatrix} &\define \scalemath{0.95}{ \st {\rm col}\left\{\begin{bmatrix}
\bar{\g}^\ell(\bzeta_{i}) \\ 
\check{\g}^\ell(\bzeta_{i})
 \end{bmatrix} \right\} = \begin{bmatrix}
{\rm col}\big\{
\bar{\g}^\ell(\bzeta_{i})
  \big\} \\ 
{\rm col}\big\{ 
\check{\g}^\ell(\bzeta_{i}) \big\}
 \end{bmatrix} } \\
 \begin{bmatrix}
\bar{\f}^\ell(\bsw_{i-1}) \\ 
\check{\f}^\ell(\bsw_{i-1})
 \end{bmatrix} &\define \scalemath{0.95}{\st {\rm col}\left\{\begin{bmatrix}
\bar{\f}^{\ell}(\bsw_{i-1}) \\ 
\check{\f}^{\ell}(\bsw_{i-1})
 \end{bmatrix} \right\} = \begin{bmatrix}
{\rm col}\big\{
\bar{\f}^{\ell}(\bsw_{i-1})
  \big\} \\ 
{\rm col}\big\{ 
\check{\f}^{\ell}(\bsw_{i-1}) \big\}
 \end{bmatrix} }\\
 \begin{bmatrix}
\bar{\v}_{i} \\ 
\check{\v}_{i}
 \end{bmatrix} &\define \st {\rm col}\left\{\begin{bmatrix}
\bar{\v}_{i}^\ell \\ 
\check{\v}_{i}^\ell
 \end{bmatrix} \right\} = \begin{bmatrix}
{\rm col}\left\{
\bar{\v}_{i}^\ell
  \right\} \\ 
{\rm col}\left\{ 
\check{\v}_{i}^\ell \right\}
 \end{bmatrix} \label{calT-noise}
}
 Since $\cJ$ is block diagonal, the operation $\st \mathcal{J}\tran \st\tran $ performs a similar reordering with respect to the diagonal blocks -- Figure \ref{fig:TXT} illustrates this operation visually.
\noindent Thus,
\eq{
\st \mathcal{J} \st\tran  
&=\begin{bmatrix}
I_M & 0 \\
0 & {\rm blkdiag}\{\check{J}_\ell \otimes I_{M_\ell}   \} \end{bmatrix} 
}
plugging into \eqref{errorT} we arrive at thew following conclusion.
\begin{lemma} {\bf (Transformed error recursion)}: Following similarity and permutation  transformations, the error recursion \eqref{error1} can be transformed into the following form 
\eq{
\begin{bmatrix}
\bw_{i} \\ 
\cw_i
 \end{bmatrix} &=   \begin{bmatrix}
\bw_{i-1}+\bar{\g}(\bzeta_{i})+\bar{\f}(\bsw_{i-1})+\bar{\v}_{i} \\ 
\check{\cJ}\big(\cw_{i-1} +\check{\g}(\bzeta_{i})+\check{\f}(\bsw_{i-1})+\check{\v}_{i} \big)
 \end{bmatrix}
\label{errorT-scaled}}
where 
\eq{ 
\check{\cJ} \define {\rm blkdiag}\{\check{J}_\ell \otimes I_{M_\ell}   \}
}
\qd
\end{lemma}
\noindent 
We are now ready to state the main result regarding the coupled diffusion algorithm \eqref{GM Stochastic Diffusion}.
\subsection{Mean-Square Convergence} 
\begin{theorem}(\textrm{\bf{Mean-square convergence)}}: Under Assumptions \ref{feasible assump}--\ref{noisemodel:assump}, the coupled diffusion algorithm \eqref{GM Stochastic Diffusion} converges in the mean-square-error sense for sufficiently small step-sizes $\mu$ (see \eqref{step-size(all)}), namely, it holds that for $i \geq 0$
\eq{
\begin{bmatrix}
\Ex \|\bw_{i}\|^2 \\ 
\Ex \|\cw_{i}\|^2
 \end{bmatrix} \preceq  \Gamma \begin{bmatrix}
\Ex \|\bw_{i-1}\|^2 \\ 
\Ex \|\cw_{i-1}\|^2
 \end{bmatrix}  +\begin{bmatrix}
c_1\\ 
c_2
 \end{bmatrix}  \label{theorem-meansquare_gamma}
}
where $\Gamma$ is a stable matrix and $\{c_1,c_2\}$ are independent of time (see \eqref{gamma_c_expressions}). It follows that, for every agent $k$,
\eq{
\limsup\limits_{i\rightarrow \infty}
\Ex \|w^{\ell,\star}- \w^\ell_{k,i} \|^2 &\leq O(\mu)+O(\mu^2 \eta^4)
\label{theorem-meansquare-agent-k-var-l}} 
for all $\ell \in \cI_k$.
\qd \label{theorem-meansquare}
\end{theorem}

Proof: See Appendix \ref{appendix-proof}.

Theorem \ref{theorem-meansquare} means that the expected squared distance between $\w_{k,i}^\ell$ and $w^{\ell,\star}$ is upper bounded by some value on the order of $\mu$ or $\mu^2 \eta^4$, whichever is larger. This implies that we can get arbitrarily close to the optimal penalized solution $w^\star={\rm col}\{w^{\ell,\star}\}_{\ell=1}^L$ by choosing $\mu$ arbitrarily small. Moreover, from Theorem \ref{theorem-approachoptimal}, we can get arbitrarily close to the original problem \eqref{glob2} by choosing $\eta$ arbitrarily large. From the step size condition \eqref{step-size(all)} we see that $\mu<O(1/\eta^2)$, therefore we can choose $\eta={c / \mu^\theta}$ for some constant $c$ and $0<\theta<0.5$. This way the problem will depend on $\mu$ only and as $\mu \rightarrow 0$, the iterates $\{\w_{k,i}^\ell\}_{k \in \cC_\ell}$ approach the optimizer of the original problem $w^{\ell,o}$ asymptotically.
Another conclusion from Theorem \ref{theorem-meansquare} is that the convergence rate is upper bounded by  the spectral radius of the matrix $\Gamma$ (from \eqref{conv_rate} and ignoring the $\eta$ terms):
\eq{
\rho(\Gamma) = \max \{1-\mu \nu +O(\mu^2),\lambda(2)+O(\mu)\}
} 
 where $\lambda(2)=\max_{ \ell \in \{1,\cdots,L\}} | \lambda_{\ell}(2)|$ and $\lambda_{\ell}(2)$ is the second largest eigenvalue in magnitude of the combination matrix $A_\ell$ (the largest eigenvalues is equal to one). The smaller $\lambda_{\ell}(2)$ is, the more connected $\cC_\ell$ is.  Apart from reducing communication and memory allocation, this result shows the importance of solving \eqref{glob2} directly and how the clusters affect the convergence rate, i.e., the convergence rate is directly affected by the connectivity of the clusters instead of the network. 
 \begin{remark}\label{remark-tighter_bounds}{ \rm 
Note that the $O(\mu)$ term in \eqref{theorem-meansquare-agent-k-var-l} is due to the persistence gradient noise component. In adaptive systems, constant step sizes are used to allow the algorithm to track changing minimizers. For example, when the distribution of the streaming data changes, the minimizer also changes and if a decaying step size is used then the algorithm will lose track of the minimizer as the step size approaches zero. This means that in practice we only need to choose a sufficiently small step size and sufficiently large penalty factor. For example, in the flow problem given in \eqref{appl-flow-penalty}, the penalty factor $\eta$ is set to be  large enough so that under constant step size and slowly varying flows $\{\b_k(i)\}$, the algorithm is still able to track the changing minimizer. In general, the penalized problem optimizer approaches the minimizer of the original problem as $\eta$ approaches infinity. Under some additional assumptions, an exact differentiable penalty function can be constructed \cite{huyer2003new} and, therefore, there exists scenarios such that the minimizers of problems \eqref{glob2} and \eqref{penalized_cost} approach each other for large enough $\eta< \infty$. 
\qd
}
\end{remark}
\begin{remark}\label{remark-th}{ \rm
 In practice, the $O(\mu^2 \eta^4)$ term in \eqref{theorem-meansquare-agent-k-var-l} is tighter and $\mu$ can be chosen to satisfy $\mu <O(1/\eta)$ instead of $\mu<O(1/\eta^2)$. Although unnecessary for the convergence analysis,
if desired, this can be tightened to
$O(\mu^2 \eta^2)$ by calling upon the following observation. The optimality condition of \eqref{penalized_cost} is:
\eq{
 0&=\sum_{k \in \mathcal{C}_\ell} \big( \grad_{w^\ell}J_{k}(w_k^\star)+\eta \grad_{w^\ell}p_{k}(w_k^\star) \big)  \nonumber \\  
 &=\sum_{k \in \mathcal{C}_\ell} \bigg( \grad_{w^\ell}J_{k}(w_k^\star) + \sum_{u=1}^{U_k} \eta \grad\delta^{{\rm EP}}(h_{k,u}(w_k^\star)) \grad_{w^\ell}h_{k,u}(w_k^\star)  \nonumber \\
 & \quad +\sum_{v=1}^{V_k} \eta \grad\delta^{{\rm IP}}(g_{k,v}(w_k^\star)) \grad_{w^\ell}g_{k,v}(w_k^\star) \bigg), \ \forall \ \ell \label{lagrangian}
} 
Assume the optimal value $w^o$ is a regular point for the constraints, meaning that the gradients of the equality constraints and the active inequality constraints $\{\grad_{w} h_{k,u}(w^o),\grad_{w} g_{k,v'}(w^o)\}$ are linearly independent (where an active constraint means that $g_{k,v'}(w_k^o)=0$ for some $v'$ where $w_k^o={\rm col}\{w^{\ell,o}\}_{\ell \in \cI_k}$). Then, it is shown in \cite[pp.~479-481]{43} (see also \cite[pp.~392-393]{bertsekas1999nonlinear}):
\eq{
\eta \grad\delta^{{\rm EP}}(h_{k,u}(w_k^\star))  \rightarrow y^o_{k,u}, \ 
 \eta \grad\delta^{{\rm IP}}(g_{k,v}(w_k^\star)) \rightarrow z^o_{k,v}
}
as $\eta \rightarrow \infty$. Here, the variables $y^o_{k,u}$ and $z^o_{k,v}$ correspond to the optimal unique dual variables of problem \eqref{glob2} associated with the constraints $h_{k,u}(.)$ and $g_{k,v}(.)$. Thus, it holds that each $\eta \grad_{w^\ell}p_{k}(w_k^\star)$ converges. We know that any convergent sequence is bounded and, hence, $\|\eta \grad_{\ssw}\cP(\sw^\star)\|$ is bounded by some constant independent of $\eta$. This observation can be used in the proof of Theorem \ref{theorem-meansquare} to tighten the bound in \eqref{theorem-meansquare-agent-k-var-l} to $O(\mu)+O(\mu^2 \eta^2)$ (see Remark \ref{remark-proof}). \qd
}
\end{remark}
\begin{figure*}[!t]
\vspace{-5mm}
\centering
\begin{subfigure}[]{0.32\textwidth}
\includegraphics[scale=0.34]{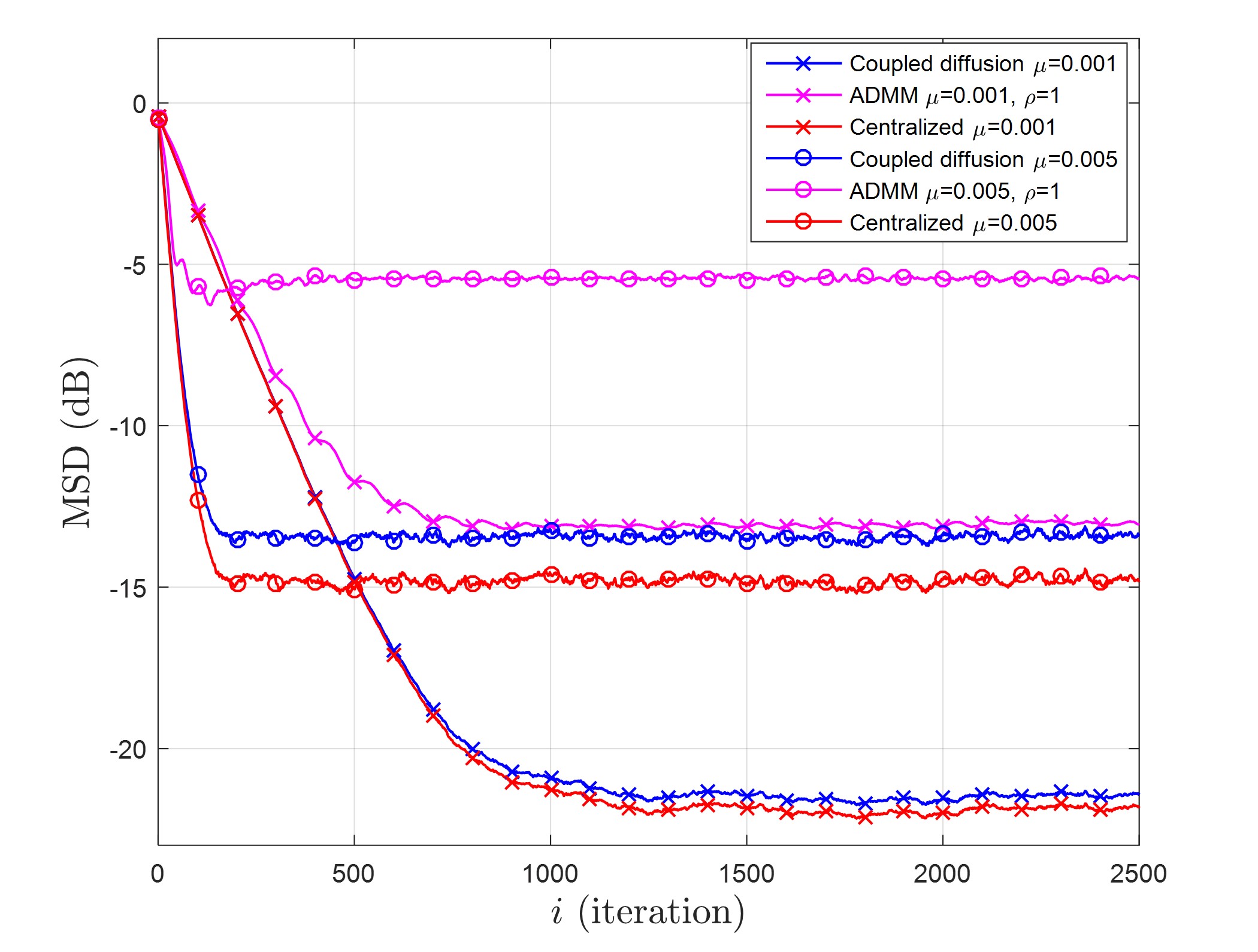}
\caption{}
\label{fig:simulation_unconstrained}
\end{subfigure}
\begin{subfigure}[]{0.32\textwidth}
\includegraphics[scale=0.34]{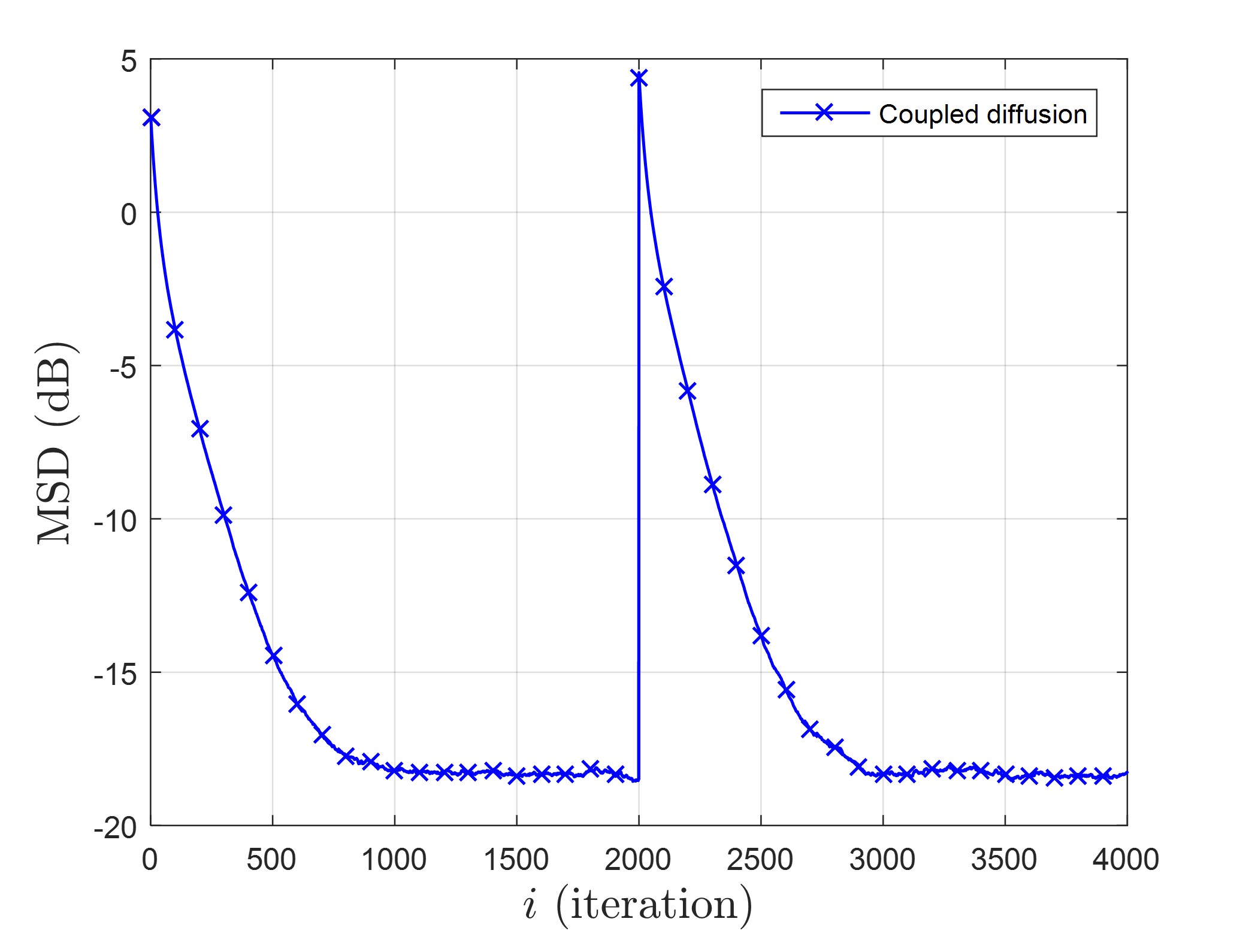}
\caption{}
\label{fig:simulation_adaptive_constrained}
\end{subfigure}
\begin{subfigure}[]{0.32\textwidth}
\includegraphics[scale=0.34]{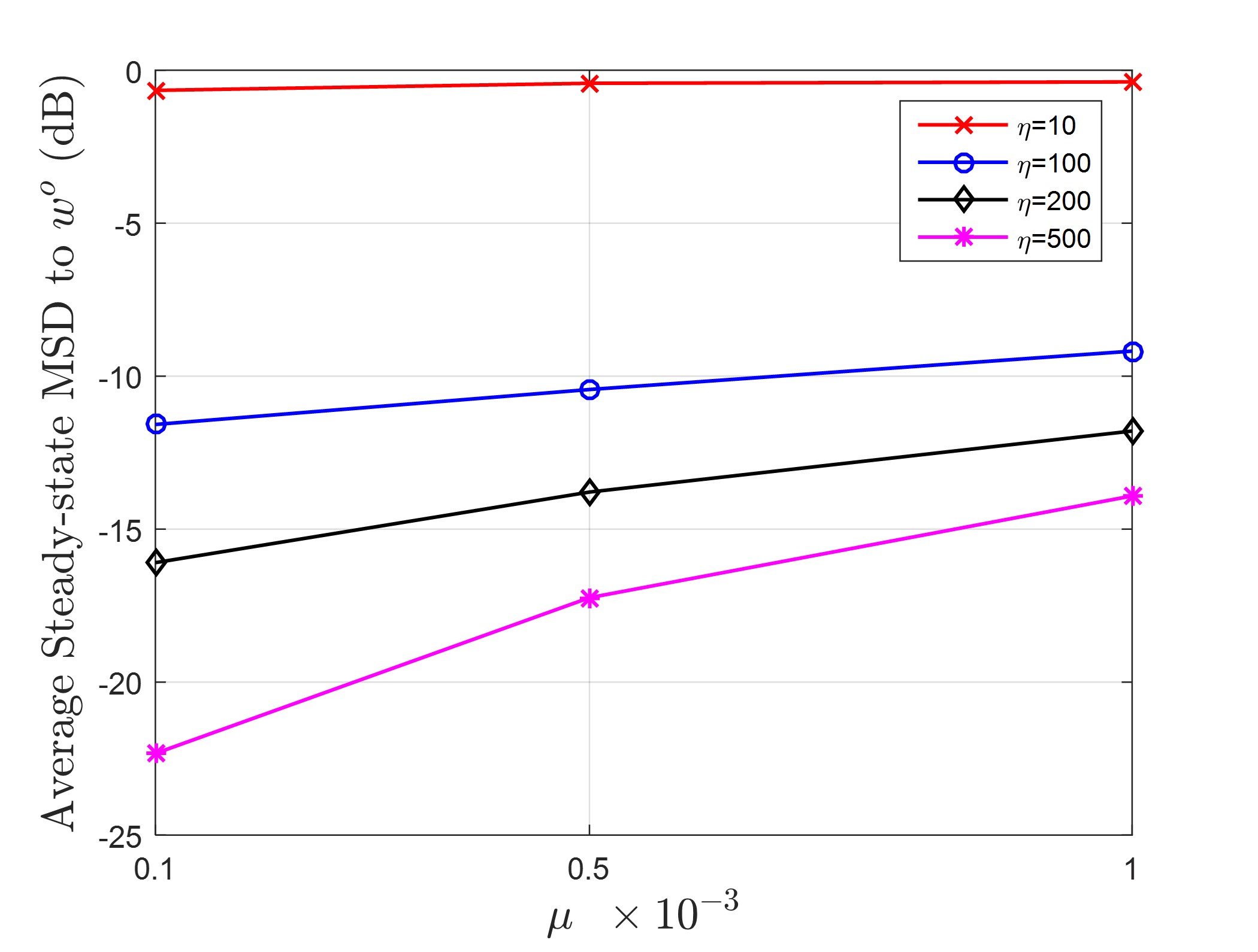}
\caption{}
\label{fig:eta_mu}
\end{subfigure}
\vspace{1mm}
\caption{ \small (a) Network MSD for two different step-sizes. (b) Network MSD for $\mu=0.001$ and $\eta=100$ showing the adaptive coupled diffusion algorithm where in iteration $2000$ the minimizer $w^\star$ changes by randomly regenerating the constraints. (c) Average steady-state MSD for different values of $\mu$ and $\eta$.}
\end{figure*}
\section{Example and Simulation Results}
In this section, we illustrate the performance of the coupled diffusion strategy \eqref{stochastic-diff(a)}--\eqref{stochastic-diff(c)} for a least-squares model fitting problem under streaming data.

Let $w={\rm col}\{w^1,\cdots,w^L\} \in \real^{M}$ and consider $N$ agents with distributed cost:
\eq{
\scalemath{0.96} {J^{{\rm glob}}(w)=\Ex \left\|\H_{i}w-\y_i\right\|^2
=\sum_{k=1}^N \Ex \big(\h_{k,i}\tran w - \y_k(i)\big)^2 }
\label{sim-glob}}
where $\h_{k,i} \in \real^{M}$ and $\y_k(i) \in \real$ are the $k$-th row of $\H_{i} \in \real^{N \times M}$ and the $k$-th element of $\y_i \in \real^{N}$, respectively. Assume $\H_{i} \in \real^{N \times M}$ and $\y_i \in \real^{N}$ are related via the linear model:
\eq{
\y_i=\H_{i}w^\bullet+ \v_i
}
with $\v_i \in \real^{N}$ representing a random Gaussian noise independent of $\H_{i}$ with covariance $\Sigma={\rm diag}\{\sigma_{v,k}(i)\}_{k=1}^N$. Assume the features $\{\h_{k,i}\}$ are sparse in the sense that they are zero at the location of $\{w^\ell\}$ if $\ell \notin \cI_k$ where $\cI_k$ represent the indices where the features are nonzero. Specifically, if we divide each feature as $\h_{k,i}=[\h_{k,i}^1,\cdots,\h_{k,i}^L]$, then we assume $\h_{k,i}^\ell=0_{M_\ell}$ if $\ell \notin \cI_k$. Thus, if we let 
\eq{
w_{k} \define{\rm col}\{w^\ell\}_{\ell \in \cI_k}, \quad \overline{\h}_{k,i}\define {\rm col}\{\h_{k,i}^\ell\}_{\ell \in \cI_k} 
}
then, we can rewrite \eqref{sim-glob} as 
\eq{
J^{{\rm glob}}(w)=\sum_{k=1}^N \Ex \left(\overline{\h}_{k,i}\tran w_k - \y_k(i)\right)^2 
\label{sim-glob2}}
Problems with a global function of the type \eqref{sim-glob2} naturally arises when different subsets of data are dispersed over $N$ processors (or agents) \cite[pp.~53]{26}. Also, in robust power system state estimation as in \eqref{power-system-estim} with streaming data instead.
\subsection{Unconstrained Case}
\begin{figure}[H]
\centering
\includegraphics[scale=0.65]{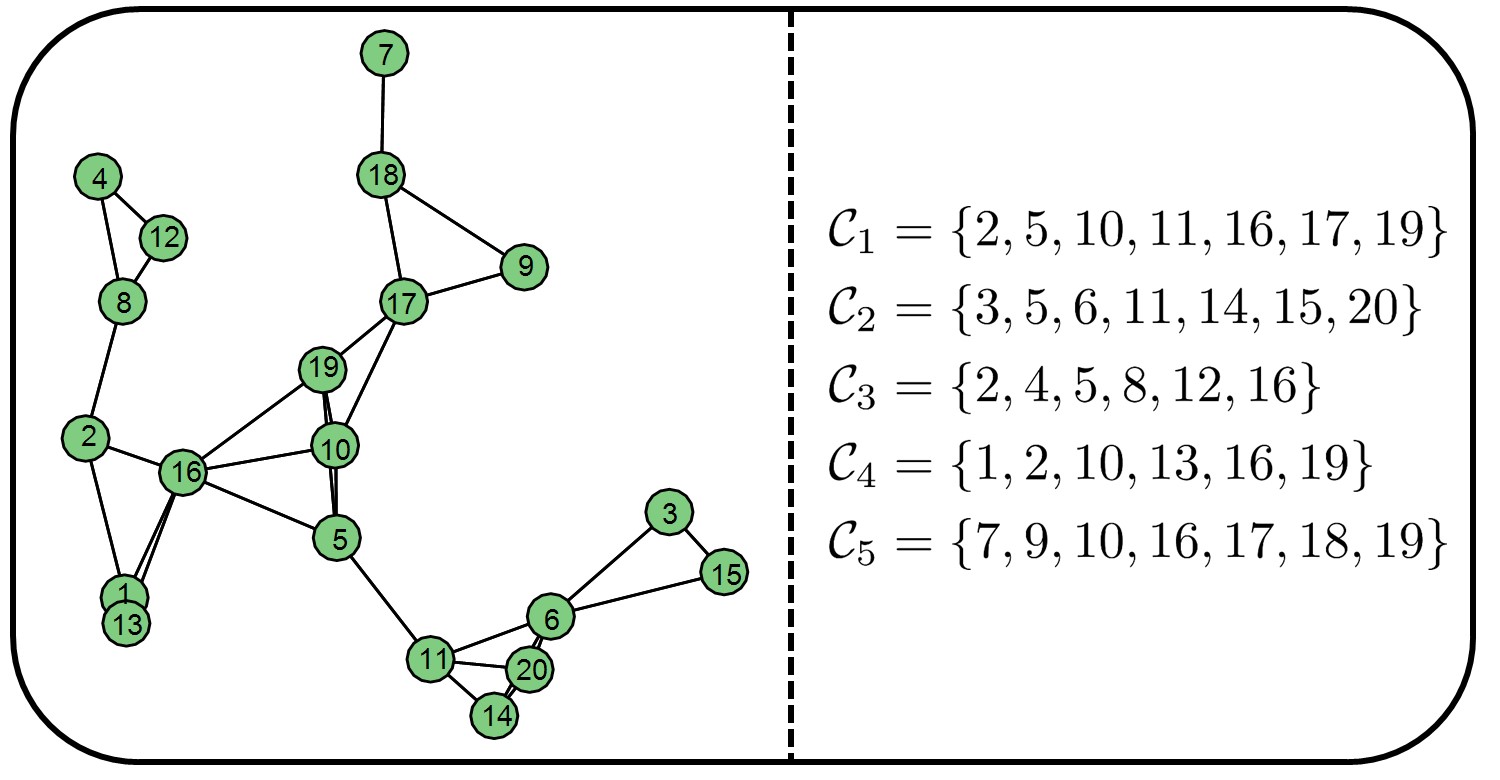}
\caption{Network topology and clusters used in the simulation.}
\label{fig:simulation_network}
\end{figure}

We first show the performance of the proposed algorithm for an unconstrained case by minimizing the cost \eqref{sim-glob2}. In our simulation we consider the network of $N=20$ agents shown in Figure \ref{fig:simulation_network}. The global vector $w$ is of size $M=25$. The number of sub-vectors is $L=5$ each of size $5 \times 1$ and each cluster is given in Figure \ref{fig:simulation_network}. The model parameter $w^\bullet$ is chosen randomly and normalized to one. The noise variances $\{\sigma_{v,k}\}$ are chosen uniformly at random between $-20$ and $-30$ dB. The covariance matrix $\Ex \overline{h}_{k,i} \overline{h}_{k,i}\tran=R_{h,k} \in \real^{Q_k \times Q_k}$ is generated as $R_{h,k}=U_k \Lambda_k U_k\tran$ where $U_k$ is a randomly generated orthogonal matrix and $\Lambda_k$ is a diagonal matrix with each diagonal entry uniformly chosen between $1$ and $3$. To compare the performance with other algorithms, we simulate a linearized version of the ADMM approach from \cite{26}. The algorithm from \cite{26} is:
\eq{
w_{k,i+1}&=\argmin_{w_k} \left(J_k(w_k)+y_{k,i}\tran w_k+{\rho \over 2}\|w_k-z_{k,i}\|^2 \right)\label{admm1} \\
z^\ell_{i+1}&={1 \over N_\ell} \sum_{k \in \cC_\ell} (w_{k,i+1}^\ell + {1 \over \rho} y_{k,i}^\ell), \ \forall \ \ell=1,\cdots,L \label{admm2}\\
y_{k,i+1}&=y_{k,i}+\rho (w_{k,i+1}-z_{k,i+1})\label{admm3}
}
where $y_k={\rm col}\{y_k^\ell\}_{\ell \in \cI_k}$, $z_{k}=\{z^\ell\}_{\ell \in \cI_k}$, and $\rho>0$. Now note that unlike the coupled diffusion strategy, this algorithm is not a first order algorithm and step \eqref{admm1} requires an inner iteration loop unless a closed form solution exists. Moreover, under adaptive networks it is not possible to solve step \eqref{admm1} using for example gradient descent with constant step-size due to the unknown cost $J_k(w_k)$, and a decaying step-size is required, which only converges asymptotically. Therefore, in our simulation for comparison we linearize this step by employing one stochastic gradient step with constant step-size $\mu$. Note also that step \eqref{admm2} requires global knowledge. We also simulate the centralized recursion \eqref{central-inc1}--\eqref{central-inc2} using stochastic gradient and with diagonal scaling $D={\rm diag}\{{1 \over N_\ell} I_{M_\ell}\}_{\ell=1}^L$ used to make the convergence rate practically the same. 

The simulation result is shown in Figure \ref{fig:simulation_unconstrained}, which plots the instantaneous network MSD
\eq{
{\rm MSD}(i)=\sum\limits_{\ell=1}^L{1\over N_\ell}\sum\limits_{k \in \cC_\ell} \Ex \|w^{\ell,\star}-\w^{\ell}_{k,i}\|^2
}
for different step sizes. The combination matrices $A_\ell$ are chosen using the Metropolis rule \eqref{Metropolis}. We see that the coupled diffusion algorithm outperforms its ADMM counterpart even though the ADMM uses global information in step \eqref{admm2}, which indicates that primal-dual methods do not necessarily perform well under adaptive networks as shown in \cite{41}. We also notice that the smaller the step-size is, the smaller is the steady-state MSD, and, moreover, the closer the coupled diffusion becomes to the centralized. 
\subsection{Constrained Case}
 We now add constraints to the previous setting where now the network is interested in solving a linearly constrained least-squares with the cost given in \eqref{sim-glob2} and $L$ linear constraints $
Gw=b
$,  
where $G \in \real^{L \times M}$ and $b \in \real^{L}$. Let $g_\ell\tran$ and $b_\ell$ denote the $\ell$-th row of $G$ and the $\ell$-th element in $b$, then we assume that only one agent $k_c$ in cluster $\cC_\ell$ is aware of the $\ell$-th constraint, which has the form:
\eq{
g_\ell\tran w = g_{\ell,k_c}\tran w_{k_c} = b_\ell
}
where $g_{\ell,k_c} \in \real^{Q_k}$. This means that $g_\ell$ has zero entries at the location of $w^\ell$ if $\ell \notin \cI_{k_c}$. This situation occurs for example when agent $k_c$ is the decision making agent regarding variable $w^\ell$ with knowledge about the constraint concerning $w^\ell$, while the other agents are just observation agents that collects data related to $w^\ell$.

In our simulation, $\{g_{\ell,k_c}\}$ are generated using the standard Gaussian distribution and normalized to one and $\{b_\ell\}$ are generated uniformly between $-1$ and $1$. The agents that have knowledge about the constraints $\ell=\{1,2,3,4,5\}$ are $k_c=\{2,10,16,5,17\}$, respectively. We used the quadratic penalty $\delta^{\rm EP}(x)= x^2$ to penalize the constrained. Figure \ref{fig:simulation_adaptive_constrained} shows the MSD for the coupled diffusion algorithm with $\mu=0.001$ and $\eta=100$. In order to illustrate that the algorithm is capable of tracking a changing minimizer, we randomly regenerate $\{g_{\ell,k_c}\}$ and $\{b_\ell\}$ at $i=2000$ so that the constraints changes and, thus, $w^\star$ changes. We see that the algorithm is capable of tracking a changing minimizer. In Figure \ref{fig:eta_mu}, we plot the steady-state MSD to the constrained problem minimizer $w^o$ for different values of step-sizes $\mu$ and penalty parameters $\eta$. It is observed that for relatively small penalty parameter $\eta=10$, the MSD becomes unaffected by how small $\mu$ is. This is because for small $\eta$, the penalized optimal vector $w^\star$ is not a good approximation of $w^o$. However, for large penalties, we notice that the smaller the step-size $\mu$ is, the smaller the MSD becomes. This is because from Theorem \ref{theorem-approachoptimal} we know that for sufficiently large $\eta$, the penalized optimal $w^\star$ becomes close to $w^o$ and from Theorem \ref{theorem-meansquare}, we know that the smaller the value of $\mu$ is, the closer the coupled diffusion iterates becomes to $w^\star$.
\section{Concluding Remarks}
In this work, we developed a distributed optimization algorithm that solves a general sum of local cost functions subject to the intersection of all local constraints, where each local cost and constraint may contain partial blocks of the global variable. We proved convergence under constant-step sizes in stochastic scenarios, where we show that we can get arbitrarily close to the optimal solution. Constant-step sizes are important in adaptation and learning where the optimizer location may drift with time. We also showed how the coupling across agents affects the convergence rate of the designed algorithm.
\begin{appendices}
 \section{Useful Bounds}  \label{appendix-bounds}
 In this appendix we give bounds that will be used in Appendix \ref{appendix-proof} for the proof of the main convergence Theorem \eqref{theorem-meansquare}.
 \begin{lemma}
Under Assumptions \ref{feasible assump}--\ref{penalty-assump} and if
\eq{
\mu < {1 \over \nu+N(\delta+\eta \delta_p) } \label{step-size(a)} 
}
then it holds that:
\eq{
\|\widetilde{w}_{i-1}+\mu \grad_w J^{\rm glob}_{\eta}(w_{i-1})\|^2 \leq (1-\mu \nu)^2 \|\widetilde{w}_{i-1}\|^2
 \label{lemma_inequal_3}}
where $\widetilde{w}_{i-1}\define w^\star-w_{i-1}$ and
\eq{
 \delta \define \max_{k} \delta_{k} , \quad
 \delta_{p} \define \max_{ k } \delta_{p,k}  
 }
\end{lemma}
\noindent {\bf Proof:}   The proof follows from \cite[~Lemma 10]{qu2017harnessing} with the Lipschitz constant set to $N(\delta+\eta \delta_p)$. We now show that the Lipschitz constant of $J^{\rm glob}_{\eta}(.)$ is upper bounded by $N(\delta+\eta \delta_p)$. Under Assumptions \ref{cost-assump}--\ref{penalty-assump} and for any $x$ and $y$ of the same structure as $w$, it holds that:
\eq{
& \big\|\grad_w J^{\rm glob}_{\eta}(x)-\grad_w J^{\rm glob}_{\eta}(y) \big\|  \nonumber \\
&  = \left\|  \sum_{k =1}^N \big(  \grad_{w}J_{k,\eta}(x_{k})- \grad_{w} J_{k,\eta}(y_{k}) \big)\right\| \nonumber \\
&  \overset{(a)}{\leq}       \sum_{k =1}^N \big\|  \grad_{w_k}J_{k,\eta}(x_k)- \grad_{w_k} J_{k,\eta}(y_k) \big\| \nonumber \\
& \leq    \sum_{k =1}^N (\delta_{k}+\eta \delta_{p,k}) \big\|  x_k- y_k \big\| \leq    N (\delta+\eta \delta_p) \big\|  x- y \big\|
}
where in step (a) we used the triangle inequality and the fact that the block entries $\{\grad_{w^\ell} J_{k,\eta}(.)\}$ are zero if $k \notin \cC_\ell$.
 \qd

 \noindent The following simple fact will be useful for the proof of the next lemma. 
 
\noindent  {\bf Fact 1.}
 {\em  Let $x_{k}={\rm col}\{z_{k}^\ell\}_{\ell \in \cI_k}$, $\ssx^\ell={\rm col}\{z_{k}^\ell\}_{k \in \cC_\ell}$ , and $\ssx={\rm col}\{\ssx^\ell\}_{\ell=1}^L$. Then, it holds that:
 \eq{
 \|\ssx\|^2=\sum_{\ell=1}^L \sum_{k \in \cC_\ell} \|z^\ell_k\|^2 =  \sum_{k=1}^N \sum_{\ell \in \cI_k} \|z^\ell_k\|^2 =  \sum_{k=1}^N \|x_k\|^2 \label{fact}
 } 
 }
{\bf Proof:} The proof follows from fact that the squared norm of a vector is equal to the sum of squared norms of its individual blocks.
 \qd
 
 For the statement and proof of the next Lemma we introduce the centroid vectors:
 \eq{
\w_{c,i} &\define{\rm col}\left\{ \w^{\ell}_{c,i} \right\}_{\ell=1}^L \in \real^{M}  \\
\w'_{k,i} &\define{\rm col}\left\{ \w^{\ell}_{c,i} \right\}_{\ell \in \cI_k} \in \real^{Q_k}
\label{average_itertes}}
where $\w^{\ell}_{c,i} \define \sum_{s \in \cC_\ell} r_\ell(s) \w_{s,i}^\ell \in \real^{M_\ell}$ is the centroid of the estimates $\{\w_{s,i}^\ell\}_{s \in \cC_\ell}$ of the variable $w^\ell$.
 \begin{lemma}
  Under Assumptions \ref{cost-assump} and \ref{penalty-assump}, the following bounds hold:
\eq{
 & \big\| \bar{\g}(\bsw_{i-1})+ \bar{\f}(\bsw_{i-1})-\mu \grad_w J^{\rm glob}_{\eta}(\w_{c,i-1})\big\|^2  \nonumber \\
  & \quad \leq \mu^2  (\delta+\eta \delta_p)^2  v_1 N_x \|\cw_{i-1} \|^2
  \label{lemma_inequal_1}}
  and
  \eq{
  &\big\|\bar{\g}(\bzeta_{i})-\bar{\g}(\bsw_{i-1})\big\|^2 \leq 2 \mu^4 \eta^2 \delta^2 N_x n  \|\grad_{\ssw} \cP(\sw^\star)\|^2 \nonumber \\
 & \quad +
2 \mu^4 \eta^2 \delta^2 \delta_p^2 v_2 N_x n ( \|\bw_{i-1}\|^2+\|\cw_{i-1}\|^2)
 \label{lemma_inequal_2}}
where $v_1 = \max_\ell (V_{2,\ell} \otimes I_{M_\ell})^2$, $v_2\define\|\cV^{-1}\|^2$, $N_x=\max_\ell\{N_\ell\}$, and $n=\|\cR^{-1}\|^2$.
 \end{lemma}
\noindent {\bf Proof:} We first show \eqref{lemma_inequal_1}. Note that from \eqref{g-bar-check}:
\eq{
\bar{\g}^\ell(\bsw_{i-1}) &= \mu (r_\ell\tran R_\ell^{-1} \otimes I_{M_\ell}) \grad_{\ssw^\ell} \cJ(\bsw_{i-1}) \nonumber \\
&= \mu(\one_{N_\ell}\tran \otimes I_{M_\ell}) \grad_{\ssw^\ell} \cJ(\bsw_{i-1}) \nonumber \\
&=\mu  \sum_{k \in \cC_\ell}   \grad_{w_k^\ell}J_{k}(\w_{k,i-1}) 
}  
Likewise for $\bar{\f}^\ell(\bsw_{i-1})$. Thus, we have:
\eq{
&\bar{\g}^\ell(\bsw_{i-1})+\bar{\f}^\ell(\bsw_{i-1}) =\mu  \sum_{k \in \cC_\ell}   \grad_{w_k^\ell}J_{k,\eta}(\w_{k,i-1})
} 
and, therefore,
\eq{
& \big\|\bar{\g}(\bsw_{i-1})+ \bar{\f}(\bsw_{i-1})-\mu \grad_w J^{\rm glob}_{\eta}(\w_{c,i-1}) \big\|^2 \nonumber \\
&  =\sum_{\ell=1}^L \left\| {\mu N_\ell \over N_\ell} \sum_{k \in \cC_\ell} \left(  \grad_{w_k^\ell}J_{k,\eta}(\w_{k,i-1})- \grad_{w^\ell} J_{k,\eta}(\w'_{k,i-1}) \right)\right\|^2 \nonumber \\
&  \leq \mu^2 N_x  \sum_{\ell=1}^L \sum_{k \in \cC_\ell} \left\|  \grad_{w_k^\ell}J_{k,\eta}(\w_{k,i-1})- \grad_{w^\ell} J_{k,\eta}(\w'_{k,i-1}) \right\|^2 \nonumber \\
&  \overset{\eqref{fact}}{=} \mu^2 N_x   \sum_{k =1}^N  \left\|  \grad_{w_k}J_{k,\eta}(\w_{k,i-1})- \grad_{w_k} J_{k,\eta}(\w'_{k,i-1}) \right\|^2 \nonumber \\
&  \overset{(a)}{\leq} \mu^2  N_x(\delta+\eta \delta_p)^2 \sum_{k =1}^N  \left\|\w_{k,i-1}- \w'_{k,i-1} \right\|^2 \nonumber \\
& \overset{\eqref{fact}}{=} \mu^2 N_x (\delta+\eta \delta_p)^2 \left\|\bsw_{i-1}-{\rm col} \{\one_{N_\ell} \otimes  \w^{\ell}_{c,i-1}\} \right\|^2 \nonumber \\
&  \overset{\eqref{clustererr=avg+dev}}{=} \mu^2  (\delta+\eta \delta_p)^2  N_x \big\|{\rm diag}\{V_{2,\ell} \otimes I_{M_\ell}\}\cw_{i-1} \big\|^2
\nonumber \\
&  \leq \mu^2  (\delta+\eta \delta_p)^2  N_x v_1 \|\cw_{i-1} \|^2
 \label{pppp_lemma00}
 }
Step (a) holds because of conditions \eqref{indv-cost-1} and \eqref{penalty-cost-1}. We now show the validity of inequality \eqref{lemma_inequal_2}. It holds that: 
\eq{
& \| \mu \eta \grad_{\ssw} \cP(\bsw_{i-1})\|^2 \nonumber \\
& = \mu^2 \eta^2  \|\grad_{\ssw} \cP(\bsw_{i-1})-\grad_{\ssw} \cP(\sw^\star)+\grad_{\ssw} \cP(\sw^\star)\|^2 \nonumber \\
& \overset{\eqref{penalty-cost-1}}{\leq} 2 \mu^2 \eta^2  \delta_p^2 \|\tsw_{i-1}\|^2+2 \mu^2 \eta^2   \|\grad_{\ssw} \cP(\sw^\star)\|^2 \label{pppp_lemma0}
}
Using an argument similar to the first five steps in \eqref{pppp_lemma00} we get:
\eq{
&\|\bar{\g}(\bzeta_{i})-\bar{\g}(\bsw_{i-1})\|^2  \leq \mu^2 \delta^2  N_x \|\bzeta_{i}-\bsw_{i-1}\|^2 \nonumber \\
& \overset{\eqref{network-recur-noise(a)}}{=} \mu^2 \delta^2 N_x  \|- \mu \eta \cR^{-1} \grad_{\ssw} \cP(\bsw_{i-1})\|^2 \nonumber \\
& \overset{\eqref{pppp_lemma0}}{\leq} 2 \mu^4 \eta^2 \delta^2  N_x n \left(\delta_p^2 \|\tsw_{i-1}\|^2+ \|\grad_{\ssw} \cP(\sw^\star)\|^2 \right) \label{pppp_lemma}
}
Note that:
\eq{
\|\tsw_{i-1}\|^2&=\big\|(\cV\tran)^{-1}\cV\tran\tsw_{i-1} \big\|^2 \leq v_2( \|\bw_{i-1}\|^2+\|\cw_{i-1}\|^2) \label{bound_wbar_wcheck}
}
Substituting the previous bound into \eqref{pppp_lemma} gives \eqref{lemma_inequal_2}.
\qd
 
 \begin{lemma} The noise terms are bounded by:
\eq{
 &\Ex  \|\bv_{i}\|^2+\Ex \|\cv_{i}\|^2   \leq   \bar{\alpha}  \left( \Ex \|\bw_{i-1}\|^2+\Ex \|\cw_{i-1}\|^2 \right) + \bar{\sigma}^2 \label{lemma_noise_bound}
}
where 
\begin{subequations} 
\eq{
\bar{\alpha} &\define \mu^2 \alpha v_2 v_3(2+4\mu^2 \eta^2  \delta_p^2 n)\\
 \bar{\sigma}^2 &\define 4 \mu^4 \eta^2  v_3 n \|\grad_{\ssw} \cP(\sw^\star)\|^2+\mu^2 v_3 \sigma^2
 } \label{alpha_sigma}
 \end{subequations}
$v_3\define\|\cV\tran \cR^{-1}\|^2$,  $\alpha \define  \max_{k}  \alpha_k$, $\alpha_k \define 2\bar{\alpha}_k$, $\sigma^2 \define \sum\limits_{k=1}^N \sigma_k^2$, and $\sigma_k^2=\bar{\sigma}_k^2+2\bar{\alpha}_k \|w_k^\star\|^2$.
\end{lemma}
\noindent {\bf Proof:}   We first note that:
\eq{ 
 & \Ex  \|\bv_{i}\|^2 +\Ex \|\cv_{i}\|^2  =\Ex \left\| \begin{bmatrix}
\bar{\v}_{i} \\ 
\check{\v}_{i}
 \end{bmatrix} \right\|^2 \overset{\eqref{calT-noise}}{=} \Ex  \left\|\st {\rm col}\left\{  \begin{bmatrix}
\bar{\v}_{i}^\ell \\ 
\check{\v}_{i}^\ell
 \end{bmatrix} \right\} \right\|^2 \nonumber \\
& =\Ex  \left\|{\rm col}\left\{\begin{bmatrix}
\bar{\v}_{i}^\ell \\ 
\check{\v}_{i}^\ell
 \end{bmatrix} \right\} \right\|^2  \overset{\eqref{v-bar-check}}{=} \Ex \big\|\mu \cV\tran \cR^{-1}   \v_i \big\|^2  \leq   \mu^2 v_3  \Ex  \| \v_i \|^2
 \label{s-bar-check}
}
 Using \eqref{noise-model(b)}, it can be easily confirmed that:
\eq{
\Ex\{\|\v_{k,i}(\bzeta_k)\|^2\mid \cf_{i} \} \leq  \alpha_k \|\tzeta_k\|^2 + \sigma_k^2 \label{noise-error-model(b)}
} 
If we take the expectation of \eqref{noise-error-model(b)} we get:
\eq{
\Ex \|\v_{k,i}(\bzeta_k)\|^2 \leq  \alpha_k \Ex\|\tzeta_k\|^2 + \sigma_k^2 \label{noise-error-model(bb)}
} 
Therefore, we can bound $\Ex \|\v_i\|^2$ as follows:
\eq{
 \Ex \|\v_i\|^2 
 &\overset{(a)}{=}    \sum_{k=1}^N \Ex  \| \v_{k,i}(\bzeta_{k,i})\|^2 
\leq   \sum_{k=1}^N ( \alpha_k \Ex \|\tzeta_{k,i}\|^2+\sigma_k^2) \nonumber \\
&\leq  \alpha  \sum_{k=1}^N \Ex \|\tzeta_{k,i}\|^2+\sigma^2  \overset{\eqref{fact}}{=} \alpha \Ex \| \tzeta_{i}\|^2 + \sigma^2  
 \label{noise-bound-zeta}}
 Step (a) holds because of \eqref{transform-vector-noise} and \eqref{fact}. Now note that:
\eq{
\Ex \| \tzeta_{i}\|^2 &\overset{\eqref{zeta-error}}{=}\Ex \big\|\tsw_{i-1} + \mu \eta \cR^{-1} \grad_{\ssw} \cP(\bsw_{i-1})\big\|^2 \nonumber \\
&\leq (2+4\mu^2 \eta^2  \delta_p^2 n) \Ex \|\tsw_{i-1}\|^2+4 \mu^2 \eta^2 n   \|\grad_{\ssw} \cP(\sw^\star)\|^2 \label{zeta-bound-w}
}
where in the last step we used Jensen's inequality and  \eqref{pppp_lemma0}. Therefore, by substituting \eqref{zeta-bound-w} into \eqref{noise-bound-zeta} and using \eqref{bound_wbar_wcheck} we conclude that
\eq{
&\Ex \|\v_i\|^2  \overset{\eqref{bound_wbar_wcheck}}{\leq}  \alpha' \left(  \Ex \|\bw_{i-1}\|^2+\Ex \|\cw_{i-1}\|^2 \right) +  \sigma'^2 \label{noise-check-bar-sum-bound}
}
where $\alpha' \define \alpha v_2(2+4\mu^2 \eta^2  \delta_p^2 n)$ and $\sigma'^2 \define 4 \mu^2 \eta^2 n   \|\grad_{\ssw} \cP(\sw^\star)\|^2+ \sigma^2$. Substituting inequality  \eqref{noise-check-bar-sum-bound} into \eqref{s-bar-check} gives the bound \eqref{lemma_noise_bound}.
 \qd 
 \section{Proof of Theorem 2}  \label{appendix-proof}
For ease of reference, we rewrite \eqref{errorT-scaled}:
\eq{
\bw_{i}&=\bw_{i-1}+\bar{\g}(\bzeta_{i})+\bar{\f}(\bsw_{i-1})+\bar{\v}_{i} \label{expanded-1}\\
\cw_{i} &= \check{\cJ}\big(\cw_{i-1} +\check{\g}(\bzeta_{i})+\check{\f}(\bsw_{i-1})+\check{\v}_{i}\big) \label{expanded-2}
}
We first bound \eqref{expanded-1}. Adding and subtracting $ \bar{\g}(\bsw_{i-1})$ from the right hand side of \eqref{expanded-1}: 
\eq{
\bw_{i}&=\bw_{i-1}+\bar{\g}_{\eta}(\bsw_{i-1})+\bar{\g}(\bzeta_{i})-\bar{\g}(\bsw_{i-1})+\bar{\v}_{i} 
}
where $\bar{\g}_{\eta}(\bsw_{i-1})\define  \bar{\g}(\bsw_{i-1})+ \bar{\f}(\bsw_{i-1})$. Similarly, adding and subtracting $\mu \grad_w J^{\rm glob}_{\eta}(\w_{c,i-1})$:
\eq{
\bw_{i}&=\bw_{i-1}+\mu  \grad_w J^{\rm glob}_{\eta}(\w_{c,i-1}) + \bar{\h}_i+\bar{\v}_{i} \label{bar_expand}
}
where we introduced
\eq{
\bar{\h}_i\define \bar{\g}_{\eta}(\bsw_{i-1})-\mu \grad_w J^{\rm glob}_{\eta}(\w_{c,i-1})+\bar{\g}(\bzeta_{i})-\bar{\g}(\bsw_{i-1})
\label{bar_h_def}}
Conditioning both sides of \eqref{bar_expand} on $\cf_{i}$, and computing the conditional second-order moments, we get:
\eq{
&\Ex \|\bw_{i}\mid \cf_{i}\|^2 = \bigg\|\bw_{i-1}+\mu \grad_w J^{\rm glob}_{\eta}(\w_{c,i-1})  +\bar{\h}_i \bigg\|^2 \nonumber \\
& \quad  + \Ex \| \bv_{i} \mid \cf_{i}\|^2     \label{bar_expand_1}
} 
where the cross term with the noise component is zero because of the gradient noise conditions given in Assumption \ref{noisemodel:assump}. Note that $\bw_{i-1}=w^\star-\w_{c,i-1}$. Appealing to Jensen's inequality we have:
\eq{
& \bigg\|\bw_{i-1}+\mu \grad_w J^{\rm glob}_{\eta}(\w_{c,i-1})  +\bar{\h}_i \bigg\|^2 \nonumber \\
 & =\bigg\|{t\over t}\bigg(\bw_{i-1}+\mu \grad_w J^{\rm glob}_{\eta}(\w_{c,i-1})\bigg)  +{1-t\over 1-t}\bar{\h}_i \bigg\|^2 \nonumber \\
 & \leq {1\over t} \bigg\|\bw_{i-1}+\mu \grad_w J^{\rm glob}_{\eta}(\w_{c,i-1}) \bigg\|^2 +{1\over (1-t)}\|\bar{\h}_i\|^2  \nonumber \\
 & \overset{\eqref{lemma_inequal_3}}{\leq} {1\over t} (1-\mu \nu)^2 \|\bw_{i-1}\|^2 +{1\over (1-t)}\|\bar{\h}_i\|^2 
}
for any $t \in (0,1)$. Note that:
\eq{
&\|\bar{\h}_i\|^2 \nonumber \\
&\overset{\eqref{bar_h_def}}{=}\big\|\bar{\g}_{\eta}(\bsw_{i-1})-\mu  \grad_w J^{\rm glob}_{\eta}(\w_{c,i-1})+\bar{\g}(\bzeta_{i})-\bar{\g}(\bsw_{i-1})\big\|^2 \nonumber \\
&\leq 2 \big\|\bar{\g}_{\eta}(\bsw_{i-1})-\mu \grad_w J^{\rm glob}_{\eta}(\w_{c,i-1})\big\|^2+2\|\bar{\g}(\bzeta_{i})-\bar{\g}(\bsw_{i-1})\|^2 \nonumber \\
&\leq 2 \mu^2 N_x \bigg(  (\delta+\eta \delta_p)^2 v_1 +2 \mu^2 \eta^2 \delta^2 \delta_p^2 v_2 n \bigg) \|\cw_{i-1} \|^2 \nonumber \\
& \quad +4 \mu^4 \eta^2 \delta^2 \delta_p^2 v_2 N_x n \|\bw_{i-1}\|^2+4 \mu^4 \eta^2 \delta^2 N_x n \|\grad_{\ssw} \cP(\sw^\star)\|^2
}
where in the last step we used \eqref{lemma_inequal_1}--\eqref{lemma_inequal_2}. Letting $t=1-\mu \nu$ and substituting the previous two bounds into \eqref{bar_expand_1}:
\eq{
&\Ex \|\bw_{i}\mid \cf_{i}\|^2 \leq \left(1-\mu \nu+{4\over   \nu} \mu^3 \eta^2 \delta^2 \delta_p^2 v_2 N_x n \right) \|\bw_{i-1}\|^2 \nonumber \\
& \quad +{2  N_x \mu \over   \nu}  \bigg(  v_1(\delta+\eta \delta_p)^2 +2 \mu^2 \eta^2 \delta^2 \delta_p^2 v_2 n \bigg) \|\cw_{i-1} \|^2 \nonumber \\
& \quad +{4  \mu^3 \eta^2 \delta^2 N_x n \over   \nu}   \|\grad_{\ssw} \cP(\sw^\star)\|^2+ \Ex \| \bv_{i} \mid \cf_{i}\|^2 \label{z-bar}
} 
We repeat a similar argument for the second relation \eqref{expanded-2}. Thus, appealing to Jensen' inequality we have:
\eq{
&\Ex \|\cw_{i} \mid \cf_{i}\|^2 \nonumber \\
&\leq \|\check{\sj}\tran \|^2 \big\|\cw_{i-1} +\check{\g}(\bzeta_{i})+\check{\f}(\bsw_{i-1})\big\|^2  +\|\check{\sj}\tran \|^2  \Ex \|\check{\v}_{i} \mid \cf_{i}\|^2 \nonumber \\
&\leq {\|\check{\sj}\tran \|^2 \over t} \|\cw_{i-1}\|^2+{ \|\check{\sj}\tran\|^2 \over 1-t}\big\|\check{\g}(\bzeta_{i})+\check{\f}(\bsw_{i-1})\big\|^2  \nonumber \\
& \quad  +\|\check{\sj}\tran \|^2  \Ex \|\check{\v}_{i} \mid \cf_{i}\|^2
\label{checkder1}}
for any arbitrary positive number $t \in (0,1)$. Now note that $\|  \check{\sj}\tran \|^2$ is equal to the spectral radius of $\check{\sj} \check{\sj}^*$ ( $(.)^*$ is the complex-conjugate transposition), and if we use the property that the spectral radius is upper bounded by any matrix norm, we have \cite[Ch. 9]{9}:
\eq{
\scalemath{0.99}{\|  \check{\sj}\tran \|^2 \leq  \|\check{\sj} \check{\sj}^*\|_1 
= \max_{\ell}\|\check{J}_{\ell} \check{J}_{\ell}^*\|_1 
=(\lambda(2)+\epsilon)^2 }
}
where  $\lambda(2)=\max_\ell | \lambda_{\ell}(2)|$, and $\lambda_{\ell}(2)$ is equal to the second largest eigenvalue in magnitude in $\check{J}_{\ell}$, which does not depend on $\epsilon$ and is strictly less than one in magnitude. From this we conclude that $\|  \check{\sj}\tran\| < 1$ for any $\epsilon < 1-\lambda(2)$.
We select $t=\|\check{\sj}\tran\| \triangleq  \rho_\epsilon $, and rewrite \eqref{checkder1} as:
\eq{
\Ex \|\cw_{i} \mid \cf_{i}\|^2 
&\leq \rho_{\epsilon} \|\cw_{i-1}\|^2+{ 1 \over 1-\rho_{\epsilon}} \big\|\check{\g}(\bzeta_{i})+\check{\f}(\bsw_{i-1}) \big\|^2  \nonumber \\
& \quad  +\rho_{\epsilon}^2  \Ex \|\check{\v}_{i} \mid \cf_{i}\|^2
\label{checkder1}}
Note that:
\eq{
&\big\|\check{\g}(\bzeta_{i})+\check{\f}(\bsw_{i-1})\big\|^2 =\sum_{\ell=1}^L \big\|\check{\g}^\ell(\bzeta_{i})+\check{\f}^\ell(\bsw_{i-1})\big\|^2 \nonumber \\
& \overset{(a)}{=}\sum_{\ell=1}^L \bigg\|(V_{1,\ell}\tran \otimes I_{M_\ell})\cR_\ell^{-1} \big(\mu \grad_{\ssw^\ell} \cJ(\bzeta_i) +\mu \eta \grad_{\ssw^\ell} \cP(\bsw_{i-1})\big)\bigg\|^2 \nonumber \\
&\leq v_4  \big\|\mu \grad \cJ(\bzeta_i) + \mu \eta \grad \cP(\bsw_{i-1})\big\|^2 \nonumber \\
&\leq 2v_4 \mu^2 \big\| \grad \cJ(\bzeta_i)\|^2 + 2v_4  \|\mu \eta \grad \cP(\bsw_{i-1})\big\|^2 \nonumber \\
&\overset{\eqref{pppp_lemma0}}{\leq} 2v_4 \mu^2 \| \grad \cJ(\bzeta_i)\|^2 \nonumber \\
& \quad + 4 v_4 \mu^2 \eta^2  \delta_p^2 \|\tsw_{i-1}\|^2+4v_4 \mu^2 \eta^2   \|\grad_{\ssw} \cP(\sw^\star)\|^2
}
 where $v_4 =\max_\ell \|(V_{1,\ell}\tran \otimes I_{M_\ell}) \cR_\ell^{-1}\|^2$ and step (a) holds from \eqref{g-bar-check}--\eqref{f-bar-check}. Note also that:
\eq{
&\| \grad \cJ(\bzeta_i)\|^2=\| \grad \cJ(\bzeta_i)-\grad \cJ(\sw^\star)+\grad \cJ(\sw^\star)\|^2 \nonumber \\
& \leq 2\| \grad \cJ(\bzeta_i)-\grad \cJ(\sw^\star)\|^2+2\|\grad \cJ(\sw^\star)\|^2 \nonumber \\
& \overset{\eqref{indv-cost-1}}{\leq} 2\delta^2\| \tzeta_i\|^2+2\|\grad \cJ(\sw^\star)\|^2
\nonumber \\
& \overset{\eqref{zeta-bound-w}}{\leq} 4\delta^2(1+2\mu^2 \eta^2  \delta_p^2 n) \Ex \|\tsw_{i-1}\|^2+8 \delta^2 \mu^2 \eta^2 n   \|\grad_{\ssw} \cP(\sw^\star)\|^2 \nonumber \\
& \quad +2\|\grad \cJ(\sw^\star)\|^2}
Substituting the last two bounds into \eqref{checkder1} and using \eqref{bound_wbar_wcheck} gives:
\eq{
& \Ex \|\cw_{i} \mid \cf_{i}\|^2 
\leq \bigg(\rho_{\epsilon}+{ 4 v_2v_4 \mu^2 \over 1-\rho_{\epsilon}}(\eta^2\delta_p^2+a )\bigg)  \|\cw_{i-1}\|^2  \nonumber \\
& \quad +{ 4 v_2v_4 \mu^2 \over 1-\rho_{\epsilon}}(\eta^2\delta_p^2+a ) \|\bw_{i-1}\|^2 + { 4v_4\mu^2   \over 1-\rho_{\epsilon}}\|\grad \cJ(\sw^\star)\|^2  \nonumber \\
& \quad + { 4 v_4 \mu^2 \eta^2 \over 1-\rho_{\epsilon}}\left(1+4\mu^2 \delta^2 n \right)    \|\grad_{\ssw} \cP(\sw^\star)\|^2  +\rho_{\epsilon}^2  \Ex \|\check{\v}_{i} \mid \cf_{i}\|^2
\label{checkder2}}
where $a\define=2\delta^2+4\mu^2 \eta^2 \delta^2 \delta_p^2 n $. If we introduce the scalar coefficients:
\begin{subequations} \label{gamma_c_expressions}
\eq{
\gamma_{11}&=1-\mu \nu+{4\over   \nu} \mu^3 \eta^2 \delta^2 \delta_p^2 v_2 N_x n+\bar{\alpha}  \\
\gamma_{12}&= {2  N_x \mu \over   \nu}  \bigg(  v_1(\delta+\eta \delta_p)^2 +2 \mu^2 \eta^2 \delta^2 \delta_p^2 v_2 n \bigg) +\bar{\alpha} \\
\gamma_{21}&=  { 4 v_2v_4 \mu^2 \over 1-\rho_{\epsilon}}(\eta^2\delta_p^2+a )+\rho_{\epsilon}^2 \bar{\alpha}    \\
\gamma_{22}&= \rho_{\epsilon}+{ 4 v_2v_4 \mu^2 \over 1-\rho_{\epsilon}}(\eta^2\delta_p^2+a )+\rho_{\epsilon}^2\bar{\alpha}  \\
c_1 &= {4   \mu^3 \eta^2 \delta^2 N_x n \over   \nu}   \|\grad_{\ssw} \cP(\sw^\star)\|^2 +\bar{\sigma}^2 \\
c_2&= { 4v_4 \mu^2 \eta^2 \over 1-\rho_{\epsilon}}\left(1+4\mu^2 \delta^2 n\right)    \|\grad_{\ssw} \cP(\sw^\star)\|^2 \nonumber \\
& \quad+{ 4v_4\mu^2   \over 1-\rho_{\epsilon}}\|\grad \cJ(\sw^\star)\|^2 +\bar{\sigma}^2
}
\end{subequations}
where $\bar{\alpha}$ and $\bar{\sigma}^2$ are defined in \eqref{alpha_sigma}. Then, by taking the expectation of \eqref{z-bar} and \eqref{checkder2} and using these parameters along with the gradient noise bound \eqref{lemma_noise_bound}, we can combine \eqref{z-bar} and \eqref{checkder2} into a single compact inequality as follows:
 \eq{
\begin{bmatrix}
\Ex \|\bw_{i}\|^2 \\ 
\Ex \|\cw_{i}\|^2
 \end{bmatrix} \preceq  \underbrace{ \begin{bmatrix}
    \gamma_{11} & \gamma_{12}   \\
\gamma_{21} & \gamma_{22}
  \end{bmatrix}}_{\define \Gamma} \begin{bmatrix}
\Ex \|\bw_{i-1}\|^2 \\ 
\Ex \|\cw_{i-1}\|^2
 \end{bmatrix}  +\begin{bmatrix}
c_1\\ 
c_2
 \end{bmatrix}  \label{gamma-recursion}
}
Note  $\Gamma$, $c_1$, and $c_2$ have entries of the following form:
\eq{
\Gamma &= \scalemath{1}{  \begin{bmatrix}
    1-O(\mu) & O(\mu)+O(\mu \eta^2)    \\
O(\mu^2 \eta^2) & \rho_\epsilon +O(\mu^2 \eta^2)
  \end{bmatrix} }
 \label{gammaorder} \\
\scalemath{1}{ \begin{bmatrix}
c_1\\ 
c_2
 \end{bmatrix}} &=  \begin{bmatrix}
O(\mu^3 \eta^2)+O(\mu^2)\\ 
O(\mu^2 \eta^2 )
 \end{bmatrix}  \label{c1c2order}
}
\begin{remark}\label{remark-proof}{\rm The property referred to in the earlier Remark \ref{remark-th} would have given $ \|\eta\grad_{\ssw} \cP(\sw^\star)\|^2=O(1)$ so that $c_1=O(\mu^2)$ and $c_2=O(\mu^2)$.}
 \qd
 \end{remark}
We continue with $c_2=O( \mu^2 \eta^2)$. Now we show that we can choose a sufficiently small $\mu$ such that $\Gamma$ is a stable matrix (i.e., the spectral radius is strictly less than one, $\rho(\Gamma) < 1$). Since the value of $\gamma_{11}$ is nonnegative under condition \eqref{step-size(a)} and the values of $\{\gamma_{12},\gamma_{21},\gamma_{22}\}$ are nonnegative, we invoke the property that the spectral radius of a matrix is upper bounded by any of its norms, and use the $1$-norm to conclude that:
\eq{
\rho(\Gamma) \leq \text{max}\{&\gamma_{11}+\gamma_{21},\gamma_{12}+\gamma_{22} \} \label{conv_rate}
}
Therefore, to find sufficient conditions that ensure the stability of $\Gamma$, the step-size $\mu$ are chosen such that:
\eq{
&1-\mu \nu+{4\over   \nu} \mu^3 \eta^2 \delta^2 \delta_p^2 v_2 N_x n + { 4 v_2v_4 \mu^2 \over 1-\rho_{\epsilon}}(\eta^2\delta_p^2+a ) \nonumber \\
& +(1+\rho_{\epsilon}^2 ) \mu^2 \alpha v_2 v_3(2+4\mu^2 \eta^2  \delta_p^2 n) < 1 \label{sumgamma1}
} 
and
\eq{
& \rho_{\epsilon}+{ 4 v_2v_4 \mu^2 \over 1-\rho_{\epsilon}}(\eta^2\delta_p^2+a )  \nonumber \\
& + {2  N_x \mu \over   \nu}  \bigg(  v_1(\delta+\eta \delta_p)^2 +2 \mu^2 \eta^2 \delta^2 \delta_p^2 v_2 n \bigg) \nonumber \\
& +(1+\rho_{\epsilon}^2 ) \mu^2 \alpha v_2 v_3(2+4\mu^2 \eta^2  \delta_p^2 n) < 1 \label{sumgamma2}
} 
 We first find conditions such that \eqref{sumgamma1} holds. Note that under condition \eqref{step-size(a)} it holds that:
 \eq{
 \mu \eta \delta_p < {1 \over N} \label{stricter}
 }
  Using this and subtracting one from both sides of \eqref{sumgamma1} we get the stricter inequality:
\eq{
&-\mu \nu+{4\over  N \nu} \mu^2 \eta \delta^2 \delta_p v_2 N_x n +{ 4 v_2v_4 \mu^2 \over 1-\rho_{\epsilon}}(\eta^2\delta_p^2+2\delta^2 a' ) \nonumber \\
& +2(1+\rho_{\epsilon}^2 )\mu^2 \alpha v_2 v_3 a' < 0 \label{sumgamma1-(b)}
}
where $a' \define 1+{2n \over N^2} $ and we used $a=2\delta^2+4\mu^2 \eta^2 \delta^2 \delta_p^2 n < 2\delta^2 a' $ due to \eqref{stricter}. Therefore, a sufficient condition to satisfy  \eqref{sumgamma1} is to choose $\mu$ such that \eqref{sumgamma1-(b)} is satisfied, which gives:
\eq{\scalemath{1}{
\mu < \frac{\nu}{ \tau_1+{4\over   N \nu}  \eta \delta^2 \delta_p v_2 N_x n   }
}\label{step-size(c)}}
where
\eq{
\tau_1 &\define { 4 v_2v_4  \over 1-\rho_{\epsilon}}(\eta^2\delta_p^2+2\delta^2 a' )+2(1+\rho_{\epsilon}^2 ) \alpha v_2 v_3 a' \label{tau1}
} 
We now find conditions such that \eqref{sumgamma2} holds. Similarly, equation \eqref{sumgamma2} can be replaced by the stricter inequality: 
\eq{
& { 4 v_2v_4 \mu^2 \over 1-\rho_{\epsilon}}(\eta^2\delta_p^2+2 \delta^2 a' )  + {2 \mu N_x \over  \nu}  \left( v_1 (\delta+\eta \delta_p)^2 +{2  \delta^2 v_2 n \over N^2}  \right) \nonumber \\
& +2(1+\rho_{\epsilon}^2 )\mu^2 \alpha v_2 v_3 a'< 1-\rho_{\epsilon} \label{sumgamma2-(b)}
} 
To simplify the expression, we let 
\eq{
\tau_2 \define {2  N_x \over  \nu}  \left( v_1 (\delta+\eta \delta_p)^2 +{2  \delta^2 v_2 n \over N^2}  \right) \label{tau2}
} 
and by using \eqref{tau1}, we can rewrite the previous inequality as:
\eq{
& \tau_2 \mu  +  \tau_1  \mu^2 < (1-\rho_\epsilon)
}
or, equivalently, by:
\eq{
- \big(\tau_2 \mu -  \tau_1 \mu^2\big)
 + 2\tau_2 \mu <(1-\rho_\epsilon) \label{232}
}
Now consider a generic inequality of the form $-(a\mu-b\mu^2)+2a\mu \leq c$. Then, we can guarantee this inequality by selecting $\mu$ to satisfy $(a\mu-b\mu^2)>0 \iff \mu < a/b$, and by also selecting $\mu$ to satisfy $2a \mu \leq c \iff \mu < c/(2a)$. Applying this conclusion to \eqref{232} we find that a sufficient condition for it to hold is to select
\eq{
\mu < \min \bigg\{\frac{\tau_2}{\tau_1 },\frac{(1-\rho_\epsilon)}{2 \tau_2} \bigg\} \label{step-size(d)}
}
  Therefore, we combine conditions \eqref{step-size(a)}, \eqref{step-size(c)}, and \eqref{step-size(d)} into the following sufficient condition:
\eq{
\mu < \scalemath{0.95}{ \min \bigg\{{1 \over \nu+N(\delta+\eta \delta_p) },\frac{\nu}{ \tau_1+{4\over   N \nu}  \eta \delta^2 \delta_p v_2 N_x n  },\frac{\tau_2}{\tau_1 },\frac{(1-\rho_\epsilon)}{2 \tau_2} \bigg\} } \label{step-size(all)}
}
Therefore, under \eqref{step-size(all)} and for sufficiently small step-sizes, it holds that (where in the denominator we are only considering the $O(\mu)$ terms):
\eq{ \limsup\limits_{i\rightarrow \infty}
 \begin{bmatrix}
\Ex \|\bw_{i}\|^2 \\ 
\Ex \|\cw_{i}\|^2
 \end{bmatrix} &\preceq   (I-\Gamma)^{-1} \begin{bmatrix}
c_1\\ 
c_2
 \end{bmatrix}
\nonumber \\
 &=   \frac{ \begin{bmatrix}
   1-\gamma_{22} & \gamma_{12}   \\
\gamma_{21} & 1-\gamma_{11}
  \end{bmatrix}}{(1-\gamma_{11})(1-\gamma_{22})-\gamma_{12}\gamma_{21}} \begin{bmatrix}
c_1\\ 
c_2
 \end{bmatrix} 
  \nonumber \\
 &= \begin{bmatrix}
O(\mu)+O(\mu^2 \eta^4 )\\ 
O(\mu^2)+O(\mu^2 \eta^2)
 \end{bmatrix}
} 
for sufficiently small step-sizes. From which we conclude that
\eq{
&\limsup\limits_{i\rightarrow \infty}
\Ex \|\tsw_i\|^2
=
\limsup\limits_{i\rightarrow \infty}
\Ex \left\|
(\cV^{-1})\tran \mathcal{T}\tran \begin{bmatrix}
 \bw_{i} \\ 
\cw_{i}
 \end{bmatrix} \right\|^2  \nonumber \\
 &\leq \big\|(\cV^{-1})\tran \big\|^2 \limsup\limits_{i\rightarrow \infty}
\bigg[\Ex \|\bw_{i}\|^2 +\Ex \|\cw_{i}\|^2\bigg] \nonumber \\
&=O(\mu)+O(\mu^2 \eta^4 )
}
for sufficiently small step-sizes. 
   \end{appendices}
\bibliographystyle{ieeetran}
\bibliography{myref,Bib_ref} 

\begin{thebibliography}{10}
\providecommand{\url}[1]{#1}
\csname url@samestyle\endcsname
\providecommand{\newblock}{\relax}
\providecommand{\bibinfo}[2]{#2}
\providecommand{\BIBentrySTDinterwordspacing}{\spaceskip=0pt\relax}
\providecommand{\BIBentryALTinterwordstretchfactor}{4}
\providecommand{\BIBentryALTinterwordspacing}{\spaceskip=\fontdimen2\font plus
\BIBentryALTinterwordstretchfactor\fontdimen3\font minus
  \fontdimen4\font\relax}
\providecommand{\BIBforeignlanguage}[2]{{%
\expandafter\ifx\csname l@#1\endcsname\relax
\typeout{** WARNING: IEEEtran.bst: No hyphenation pattern has been}%
\typeout{** loaded for the language `#1'. Using the pattern for}%
\typeout{** the default language instead.}%
\else
\language=\csname l@#1\endcsname
\fi
#2}}
\providecommand{\BIBdecl}{\relax}
\BIBdecl

\bibitem{alghunaim2018icassp}
S.~A. Alghunaim and A.~H. Sayed, ``Distributed coupled learning over adaptive
  networks,'' in \emph{Proc. IEEE ICASSP}, Calgary, Canada, April 2018, pp.
  6353--6357.

\bibitem{6}
D.~P. Bertsekas, ``A new class of incremental gradient methods for least
  squares problem,'' \emph{SIAM J. Optim.}, vol.~7, no.~4, pp. 913--926, Nov.
  1997.

\bibitem{34}
A.~Nedic and D.~P. Bertsekas, ``Incremental subgradient methods for
  nondifferentiable optimization,'' \emph{SIAM J. Optim.}, vol.~12, no.~1, pp.
  109--138, Jul. 2001.

\bibitem{35}
D.~Blatt, A.~O. Hero, and H.~Gauchman, ``A convergent incremental gradient
  method with constant step size,'' \emph{SIAM J. Optim.}, vol.~18, no.~1, pp.
  29--51, Feb. 2007.

\bibitem{36}
C.~G. Lopes and A.~H. Sayed, ``Incremental adaptive strategies over distributed
  networks,'' \emph{IEEE Trans. Signal Process.}, vol.~55, no.~8, pp.
  4064--4077, Aug. 2007.

\bibitem{56}
K.~Yuan, B.~Ying, X.~Zhao, and A.~H. Sayed, ``Exact diffusion strategy for
  optimization by networked agents,'' in \emph{{ Proc. EUSIPCO}}, Kos, Greece,
  Aug.--Sep. 2017, pp. 141--145.

\bibitem{23}
S.~S. Ram, A.~Nedic, and V.~V. Veeravalli, ``Distributed stochastic subgradient
  projection algorithms for convex optimization,'' \emph{J. {O}ptim. {T}heory
  {A}ppl.}, vol. 147, no.~3, pp. 516--545, 2010.

\bibitem{24}
A.~G. Dimakis, S.~Kar, J.~M.~F. Moura, M.~G. Rabbat, and A.~Scaglione, ``Gossip
  algorithms for distributed signal processing,'' \emph{Proc. {IEEE}}, vol.~98,
  no.~11, pp. 1847--1864, Nov. 2010.

\bibitem{25}
S.~Kar, J.~M.~F. Moura, and K.~Ramanan, ``Nonlinear observation models and
  imperfect communication,'' \emph{{IEEE} Transactions on Information Theory},
  vol.~58, no.~6, pp. 3575--3605, Jun. 2012.

\bibitem{29}
P.~Braca, S.~Marano, and V.~Matta, ``Enforcing consensus while monitering the
  environment in wireless sensor networks,'' \emph{IEEE Trans. Signal
  Process.}, vol.~56, no.~7, pp. 3375--3380, July 2008.

\bibitem{33}
S.~Kar and J.~M.~F. Moura, ``Distributed consensus algorithms in sensor
  networks: link failures and channel noise,'' \emph{IEEE Trans. Signal
  Process.}, vol.~57, no.~1, pp. 355--369, Jan. 2009.

\bibitem{5}
Z.~J. Towfic and A.~H. Sayed, ``Adaptive penalty-based distributed stochastic
  convex optimization,'' \emph{IEEE Trans. Signal Process}, vol.~62, no.~15,
  pp. 3924--3938, August 2014.

\bibitem{7}
J.~Chen and A.~H. Sayed, ``Distributed pareto optimization via diffusion
  strategies,'' \emph{IEEE J. Sel. Topics Signal Process.}, vol.~7, no.~2, pp.
  205--220, April 2013.

\bibitem{28}
------, ``On the learning behavior of adaptive networks - {P}art {I}:
  Performance analysis,'' \emph{IEEE {T}ransactions on Information Theory},
  vol.~61, no.~6, pp. 3487--3517, June 2015.

\bibitem{8}
------, ``On the learning behavior of adaptive networks - {P}art {II}:
  Performance analysis,'' \emph{IEEE {T}ransactions on Information Theory},
  vol.~61, no.~6, pp. 3518--3548, June 2015.

\bibitem{27}
A.~H. Sayed, ``Adaptive networks,'' \emph{Proceedings of the IEEE}, vol. 102,
  no.~4, pp. 460--497, Apr. 2014.

\bibitem{26}
S.~Boyd, N.~Parikh, E.~Chu, B.~Peleato, and J.~Eckstein, ``Distributed
  optimization and statistical learning via alternating direction method of
  multipliers,'' \emph{Found. Trends Mach. Lear.}, vol.~3, no.~1, pp. 1--122,
  Jan. 2011.

\bibitem{37}
J.~F. Mota, J.~M. Xavier, P.~M. Aguiar, and M.~Puschel, ``{D-ADMM}: {A}
  communication-efficient distributed algorithm for sperable optimizatoin,''
  \emph{IEEE Trans. Signal Process.}, vol.~61, no.~10, pp. 2718--2723, 2013.

\bibitem{38}
W.~Shi, Q.~Ling, K.~Yuan, G.~Wu, and W.~Yin, ``On the convergence of the {ADMM}
  in decentralized consensus optimization,'' \emph{IEEE Trans. Signal
  Process.}, vol.~62, no.~7, pp. 1750--1761, 2014.

\bibitem{39}
A.~Mokhtari, W.~Shi, Q.~Ling, and A.~Ribeiro, ``{DQM}: {D}ecentralized
  quadratically approximated alternating direction method of multipliers,''
  \emph{IEEE Trans. Signal Process.}, vol.~64, no.~19, pp. 5158--5173, 2016.

\bibitem{14}
J.~Chen, L.~Tang, J.~Liu, and J.~Ye, ``A convex formulation for learning shared
  structures from multiple tasks,'' in \emph{Proc. {ICML}}, Montreal, QC,
  Canada, June 2009, pp. 137--144.

\bibitem{15}
O.~Chapelle, P.~Shivaswmy, K.~Q. Vadrevu, S.~Weinberger, Y.~Zhang, and
  B.~Tseng, ``Multi-task learning for boosting with applications to web search
  ranking,'' in \emph{Proc. {ACM SIGKDD}}, Washington, DC, USA, Jul. 2010, pp.
  1189--1198.

\bibitem{16}
J.~Zhou, L.~Yuan, J.~Liu, and J.~Ye, ``A multi-task learning formulation for
  predicting disease progression,'' in \emph{Proc. {ACM SIGKDD}}, San Diego,
  CA, USA, Aug., 2011, pp. 814--822.

\bibitem{17}
N.~Keshava and J.~F. Mustard, ``Spectral unmixing,'' \emph{IEEE Signal Process.
  Mag.}, vol.~19, no.~1, pp. 44--57, Jan. 2002.

\bibitem{18}
J.~M. Bioucas-Dias, A.~Plaza, N.~Dobigeon, M.~Parente, Q.~Du, P.~Gader, and
  J.~Chanussot, ``Hyperspectral unimixing overview:geometrical, statistical,
  and sparse regression-based approaches,'' \emph{IEEE J. Sel. Topics Appl.
  Earth Observ.}, vol.~5, no.~2, pp. 354--379, Apr. 2012.

\bibitem{19}
R.~K. Ahuja, T.~L. Magnanti, and J.~B. Orlin, \emph{Network {F}lows: {T}heory,
  {A}lgorithms, and {A}pplications}, Prentice Hall, NJ, 1993.

\bibitem{11}
R.~Halvgaard, L.~Vandenberghe, N.~K. Poulsen, H.~Madsen, and J.~B. J{o}rgensen,
  ``Distributed model predictive control for smart energy systems,'' \emph{IEEE
  {T}rans. {S}mart {G}rid}, vol.~7, no.~3, pp. 1675--1682, April 2016.

\bibitem{20}
D.~P. Bertsekas and J.~N. Tsitsiklis, \emph{Parallel {and Distributed
  Computation: Numerical Methods}}, Prentice Hall,NJ, 1989.

\bibitem{22}
J.~Plata-Chaves, A.~Bertrand, and M.~Moonen, ``Incremental multiple error
  filtered-{X LMS} for node-specific active noise control over wireless
  acoustic sensor networks,'' \emph{in IEEE Sensor Array and Multichannel
  Signal Processing Workshop}, pp. 1--5, July 2016, Rio de Janeiro, Brazil.

\bibitem{31}
F.~Cattivelli and A.~H. Sayed, ``Distributed nonlinear {K}alman filtering with
  applications to wireless localization,'' in \emph{Proc. IEEE ICASSP}, Dallas,
  TX, Mar. 2010, pp. 3522--3525.

\bibitem{21}
V.~Kekatos and G.~B. Giannakis, ``Distributed robust power system state
  estimation,'' \emph{IEEE Trans. Power Syst.}, vol.~28, no.~2, pp. 1617--1626,
  May 2013.

\bibitem{54}
S.~A. Alghunaim, K.~Yuan, and A.~H. Sayed, ``Decentralized exact coupled
  optimization,'' in \emph{Proc. Allerton Conference on Communication, Control,
  and Computing}, Allerton, IL, October 2017, pp. 338-- 345.

\bibitem{52}
P.~D. Christofides, R.~Scattolini, D.~M. de~la Pena, and J.~Liu, ``Distributed
  model predictive control: A tutorial review and future research directions,''
  \emph{Computers and Chemical Engineering}, vol.~51, pp. 21--41, April 2013.

\bibitem{1}
R.~Nassif, C.~Richard, A.~Ferrari, and A.~H. Sayed, ``Diffusion lms for
  multitask problems with local linear equality constraints,'' \emph{IEEE
  Trans. Signal Process}, vol.~65, no.~19, pp. 4979 -- 4993, 2017.

\bibitem{12}
J.~Mota, J.~Xavier, P.~Aguiar, and M.~Puschel, ``Distributed optimization with
  local domains: {A}pplication in {MPC} and network flows,'' \emph{IEEE
  {T}rans. {A}utom. {C}ontr.}, vol.~60, no.~7, pp. 2004--2009, July 2015.

\bibitem{10}
T.~H. Summers and J.~Lygeros, ``Distributed model predictive consensus via the
  alternating directoin method of multipliers,'' in \emph{Proceedings of the
  50th Annual Allerton Confrence on Communication, Control, and Computing},
  Monticello, IL, USA, 2012, pp. 79--84.

\bibitem{32}
Y.~Pu, M.~N. Zeilinger, and C.~N. Jones, ``Inexact fast alternating
  minimization algorithm for distributed model predictive control,'' in
  \emph{IEEE Conference on Decision and Control (CDC)}, 2014, pp. 5915--5921.

\bibitem{48}
M.~Schmidt, N.~L. Roux, and F.~Bach, ``{Convergence rates of inexact
  proximal-gradient methods for convex optimization},'' \emph{in {\em Proc.}
  Conference on Neural Information Processing {\em Systems} ({NIPS})}, pp.
  6819--6824, Granada, SPAIN, 2011.

\bibitem{4}
J.~Chen, C.~Richard, and A.~H. Sayed, ``Multitask diffusion adaptation over
  networks,'' \emph{IEEE Trans. Signal Process}, vol.~62, no.~16, pp.
  4129--4144, August 2014.

\bibitem{2}
R.~Nassif, C.~Richard, A.~Ferrari, and A.~H. Sayed, ``Distributed learning over
  multitask networks with linearly related tasks,'' in \emph{Proc. of Asilomar
  Conf. on Signals, Systems, and Computers}, Pacific Grove, CA, USA, Nov. 2016,
  pp. 1390 -- 1394.

\bibitem{45}
K.~Yuan, B.~Ying, X.~Zhao, and A.~H. Sayed, ``Exact diffusion for distributed
  optimization and learning-{P}art {I}: {A}lgorithm development,'' \emph{IEEE
  Transactions on Signal Processing}, vol.~67, no.~3, pp. 708--723, Feb. 2019.

\bibitem{55}
------, ``Exact diffusion for distributed optimization and learning-{P}art
  {II}: {C}onvergence analysis,'' \emph{IEEE Transactions on Signal
  Processing}, vol.~67, no.~3, pp. 724--739, Feb. 2019.

\bibitem{arablouei2014distributed}
R.~Arablouei, S.~Werner, Y.-F. Huang, and K.~Dougancay, ``Distributed least
  mean-square estimation with partial diffusion,'' \emph{IEEE Transactions on
  Signal Processing}, vol.~62, no.~2, pp. 472--484, 2014.

\bibitem{notarnicola2017distributed}
I.~Notarnicola, Y.~Sun, G.~Scutari, and G.~Notarstefano, ``Distributed big-data
  optimization via block communications,'' in \emph{IEEE 7th International
  Workshop on Computational Advances in Multi-Sensor Adaptive Processing
  (CAMSAP)}, Curacao, Dutch Antilles, Dec. 2017, pp. 1--5.

\bibitem{47}
C.~Richard, J.~Chen, S.~K. Ting, and A.~H. Sayed, ``Group diffusion {LMS},'' in
  \emph{Proc. IEEE ICASSP}, Shanghai, China, Mar. 2016, pp. 4925--4929.

\bibitem{40}
S.~Y. Tu and A.~H. Sayed, ``Diffusion strategies outperform consensus
  strategies for distributed estimation over adaptive networks,'' \emph{IEEE
  Trans. Signal Process.}, vol.~60, no.~12, pp. 6217--6234, 2012.

\bibitem{41}
Z.~J. Towfic and A.~H. Sayed, ``Stability and performance limits of adaptive
  primal-dual networks,'' \emph{IEEE Trans. Signal Process.}, vol.~63, no.~11,
  pp. 2888--2903, 2015.

\bibitem{44}
Y.~Abbasi-Yadkori, P.~Bartlett, and A.~Malek, ``Linear programming for
  large-scale {M}arkov decision problems,'' in \emph{Proc. {ICML}}, Beijing,
  China, June 2014, pp. 496--504.

\bibitem{30}
S.~Boyd and L.~Vandenberghe, \emph{{Convex Optimization}}, Cambridge University
  Press, 2004.

\bibitem{46}
C.~K. Yu, M.~van~der Schar, and A.~H. Sayed, ``Distributed learning for
  stochastic generalized {N}ash equilibrium problems,'' \emph{IEEE Trans.
  Signal Process.}, vol.~65, no.~15, pp. 3893--3908, August 2017.

\bibitem{9}
A.~H. Sayed, ``Adaptation, learning, and optimization over neworks.''
  \emph{Foundations and Trends in Machine Learning}, vol.~7, no. 4-5, pp.
  311--801, 2014.

\bibitem{42}
B.~T. Polyak, \emph{Introduction to optimization}, NY, 1987.

\bibitem{43}
M.~S. Bazaraa, H.~D. Sherali, and C.~M. Shetty, \emph{{Nonlinear Programming:
  Theory and Algorithms}}, Wiley, NY, 1993.

\bibitem{huyer2003new}
W.~Huyer and A.~Neumaier, ``A new exact penalty function,'' \emph{SIAM Journal
  on Optimization}, vol.~13, no.~4, pp. 1141--1158, 2003.

\bibitem{bertsekas1999nonlinear}
D.~P. Bertsekas, \emph{Nonlinear programming}.\hskip 1em plus 0.5em minus
  0.4em\relax Athena scientific, Belmont, 1999.

\bibitem{qu2017harnessing}
G.~Qu and N.~Li, ``Harnessing smoothness to accelerate distributed
  optimization,'' \emph{IEEE Transactions on Control of Network Systems},
  vol.~5, no.~3, pp. 1245--1260, Sept. 2018.

\end{thebibliography}
\end{document}